\documentclass[11pt]{article}
\usepackage{amsmath}
\usepackage{cases}
\usepackage{mathrsfs}
\usepackage{amsfonts}
\usepackage{amsmath,amsthm,amssymb}

\theoremstyle{plain}
\newtheorem{theorem}{Theorem}[section]
\newtheorem{definition}{Definition}[section]
\newtheorem{lemma}{Lemma}[section]
\newtheorem{remark}{Remark}[section]
\newtheorem{proposition}{Proposition}[section]
\newtheorem{corollary}{Corollary}[section]
\newif \ifLastSection \LastSectionfalse

\numberwithin{equation}{section}

\begin{document}

	\title{\bf Compressible Navier-Stokes equations without heat conduction in $L^p$-framework}
	
	\vskip 2.5cm
\author{
 Juanzi Cai,\ \ Zhigang Wu\thanks{E-mail: zgwu@dhu.edu.cn},\ \ Mengqian Liu,\ \ \ \
\vspace{2mm}\\
\textit{\small  Department of Mathematics, Donghua University,}\\\textit{\small Shanghai 201620, P.R. China} \\
 }
\date{}
\maketitle

\textbf{{\bf Abstract:}} In this paper, we mainly consider global well-posedness and long time behavior of compressible  Navier-Stokes equations without heat conduction in $L^p$-framework. This is a generalization of Peng and Zhai \cite{peng}(SIMA, 55(2023), no.2, 1439-1463), where they obtained the corresponding result in $L^2$-framework. Based on the key observation that we can release the regularity of non-dissipative entropy $S$ in high frequency in \cite{peng}, we ultimately achieve the desired $L^p$ estimate in the high frequency via complicated calculations on the nonlinear terms. In addition, we get the $L^p$-decay rate of the solution.


\bigbreak  {\textbf{AMS Subject Classification :} 35Q30, 74H40, 76N10.
	
	
	\section{Introduction}
	
\quad\quad In this paper, our concern is
the non-isentropic compressible Navier-Stokes equations:
\begin{equation}\label{1.1}
	\left\{\begin{array}{l}
		\rho_t+{\rm div} (\rho \mathbf{v})=0, \\
		\rho[\mathbf{v}_t+(\mathbf{v}\cdot \nabla \mathbf{v})]+\nabla p(\rho,\theta)=\mu\Delta \mathbf{v}+(\mu+\lambda)\nabla{\rm div}\mathbf{v},\\
		\rho C_v[\theta_t+(\mathbf{v}\cdot\nabla)\theta]+p{\rm div}\mathbf{v}=\kappa\Delta\theta+2\mu|D(\mathbf{v})|^2+\nu({\rm div}\mathbf{v})^2,
	\end{array}
	\right.
\end{equation}
for $(t,x)\in\mathbb{R}_+\times\mathbb{R}^n(n\geq2)$. The unknowns $\rho, \mathbf{v},\theta$ and $p=p(\rho,\theta)$ represent the density, velocity, temperature and pressure of the fluid, respectively. The deformation tensor $D(\mathbf{v})=\frac{1}{2}(\nabla \mathbf{v}+(\nabla \mathbf{v})^T)$. We study the case when the coefficient of heat conduction $\kappa=0$ and the gas is ideal and polytropic, i.e., $p=\rho\theta$ and the internal energy $e=C_v\theta$. Without loss of generality, the viscosity coefficients $\mu$ and $\lambda$ are constants with $\mu>0$ and $n\lambda+2\mu>0$, and the specific heat $C_v>0$ is a constant.

Due to the significance of the physical background, this model has been widely studied in past decades. When the effect of the temperature is neglected, Eq. (\ref{1.1}) reduces to isentropic compressible Navier-Stokes equations,  and the existence, uniqueness, regularity and asymptotic behavior of its solutions can be seen in \cite{chen2,chen3,chen4,danchin1,danchin5,duan2,duan3,feireisl1,feireisl2,feireisl4,xin3,zhai,zhang} and the references therein.

When the heat conducting effect is involved, it is much more complicated due to stronger nonlinearity. Some important results have also been made in the existence, regularity and asymptotic behavior of the solutions over the past decades when $\kappa> 0$. Matsumura and Nishida \cite{matsumura} proved the existence and stability
of global small solutions around a steady-state solution.  Jiang \cite{jiang1,jiang2} considered the existence and asymptotic behavior of weak solutions for one-dimensional case and a spherically symmetric structure in dimension two and three, Hoff-Jenssen \cite{hoff} considered spherically and cylindrically
symmetric structures in the three-dimensional space. Feireisl [20] discussed the problem in
bounded domains. We also mention blowup criterions given
by Huang $et\ al.$ \cite{huang1,huang2,huang3} and Wen-Zhu \cite{wen1} and references therein, and blow-up of smooth highly decreasing at infinity solutions in Rozanova \cite{rozanova}.
There are also some important works in the framework of critical Besov spaces, see Danchin and his coauthors \cite{chikami,danchin2,danchin3,danchin4} for the local well-posedness and global well-posedness, Chen $et\ al.$ \cite{chen3} for the ill-posedness.  On the other hand, some important works were obtained in the presence of vacuum, we would like to refer to Feireisl \cite{feireisl3} for existence of variational solutions when the temperature equation is
satisfied only as an inequality in the sense of distributions, Bresch-Desjardins \cite{bresch} for the existence of global weak solutions in different conditions from \cite{feireisl3}, Wen-Zhu \cite{wen2} for the global existence of spherically and
cylindrically symmetric classical solutions to the three-dimensional case,  Huang-Li \cite{huang1} for the global existence and uniqueness of classical solutions $\mathbb{R}^3$ with small initial basic energy but
possibly large oscillations for the case that with non-vacuum far field, Wen-Zhu \cite{wen3} for the the global existence and large-time behavior of the solutions with vacuum far field if the initial mass is small.

In contrast to $\kappa>0$, there are few results in the case of $\kappa=0$ in (\ref{1.1}). It was verified that there is no global smooth solution to the Cauchy problem of the full compressible Navier-Stokes equations when the initial
density is compactly supported \cite{xin1}. Later on, Xin-Yan \cite{xin2} improved the blow-up results in \cite{xin1}. There are also some efforts on the global well-posedness of solutions. For one dimensional case,
Liu-Zeng \cite{liu} considered the system in the Lagrangian coordinates and obtained pointwise estimates and large-time
behavior of solutions by studying Green's function, Li \cite{li} discussed the global well-posedness of strong solutions in
the presence of vacuum. Comparing with the one-dimensional case, the related results for the multidimensional cases are few and some fundamental questions are still challenging. Duan-Ma \cite{duan1} investigated the global existence and convergence rates of solutions to the Cauchy problem of (\ref{1.1}) in $\mathbb{R}^3$ under the additional condition on the $L^1$-norm boundedness of the initial perturbation. Recently, this additional condition on the $L^1$-norm boundedness was removed by Chen $et al.$ in \cite{chen1}.

However, how to construct the global small solution to the Cauchy problem of (\ref{1.1}) in dimension two was open. The main reason is that $\|\nabla \mathbf{v}\|_{L^2}$ is not integrable with respect to time $t$ in two dimensional whole space, while it is crucial in \cite{duan1} for three dimensional case. Wu \cite{wugc} considered the global existence in two dimensional bounded domain due to the exponential decay rate of $\nabla \mathbf{v}$.  Recently, Peng-Zhai \cite{peng} gave a positive answer for the Cauchy problem of (\ref{1.1}) for both two and three dimensional cases in the framework of Besov space. Here they need not use the decay of $\nabla \mathbf{v}$ to get the a $priori$ estimate for non-dissipative entropy $s$ (or the temperature $\theta$). Moreover, they gave the decay rate of the solution under additional assumption of initial perturbation in negative Besov space.

To overcome the difficulties arising from the non-dissipation on $\theta$, they  in \cite{peng} rewrote system (\ref{1.1}) by selecting the new variables the pressure $p$ and entropy $s$. Based on the state equation, one has
\begin{equation}\label{1.3}
	\rho=A^{-\frac{C_v}{C_v+R}}p^{-\frac{C_v}{C_v+R}}e^{-\frac{C_v}{C_v+R}}.
\end{equation}
Then, the new system on the variables $p,\mathbf{v},s$ reads
\begin{equation}\label{1.4}
	\left\{\begin{array}{l}
		\partial_tp+\gamma p{\rm div}\mathbf{v}+\mathbf{v}\cdot\nabla p=\Gamma(\mathbf{v}), \\
		\rho\partial_t \mathbf{v}+\rho \mathbf{v}\cdot\nabla \mathbf{v}-\mu\mathbf{ v}-(\lambda+\mu)\nabla{\rm div}\mathbf{v}+\nabla p=0,\\
		\partial_ts+\mathbf{v}\cdot\nabla s=\frac{\Gamma(\mathbf{v})}{p},
	\end{array}
	\right.
\end{equation}
where $\Gamma(\mathbf{v})\triangleq 2\mu|D(\mathbf{v})|^2+\lambda({\rm div}\mathbf{v})^2$. We also consider the initial value problem of (\ref{1.4}) with the initial data
\begin{equation}\label{1.5}
	(p_0(x),\mathbf{v}_0(x),s_0(x))\rightarrow(\bar p,0,\bar s)\ \ {\rm for}\ |x|\rightarrow\infty,
\end{equation}
where the given constants $\bar p>0$ and $\bar s>0$. Let
\begin{equation*}
	\alpha_1\triangleq\sqrt{\frac{C_v}{(C_v+R)\bar\rho\bar p}},\ \  \alpha_2\triangleq\sqrt{\frac{(C_v+R)\bar p}{C_v\bar\rho}},\ \ \alpha_3\triangleq\frac{\mu}{\bar\rho},\ \
	\alpha_4\triangleq\frac{\lambda+\mu}{\bar\rho}
\end{equation*}
with $\bar\rho=p(\bar p,\bar s)$. By denoting
\begin{equation*}
	P\triangleq p-\bar p,\ \ \textbf{u}\triangleq\frac{\textbf{v}}{\alpha_1},\ \ \kappa(\rho)\triangleq\frac{1}{\rho}-\frac{1}{\bar\rho},\ \ I(P)\triangleq\frac{P}{P+\bar p},
\end{equation*}
the Cauchy problem (\ref{1.4})-(\ref{1.5}) can be rewritten as
\begin{equation}\label{1.5(0)}
	\left\{\!\!\!\begin{array}{l}
		\partial_tP+\alpha_2{\rm div}\mathbf{u}=-\alpha_1\mathbf{u}\cdot\nabla P-\frac{C_v+R}{C_v}\alpha_1P{\rm div}\mathbf{u}+\frac{C_v+R}{C_v^2}\Gamma(\alpha_1\mathbf{u}), \\
		\rho\partial_t \mathbf{u}-\alpha_3\Delta \mathbf{u}-\alpha_4\nabla{\rm div}\mathbf{u}+\alpha_2\nabla P\\
		\ \ \ \ \ \ \ \ \ \ \ \ \ \ \ \ \ \ \ \ \ \ \ \ \ \ \ \ \ =-\alpha_1\mathbf{u}\cdot\nabla \mathbf{u}-\frac{1}{\alpha_1}k(\rho)\nabla P+k(\rho)(\mu\Delta \mathbf{u}+(\lambda+\mu)\nabla{\rm div}\mathbf{u}),\\
		\partial_tS+\mathbf{u}\cdot\nabla S=\frac{R}{\bar{p}}(1-I(P))\Gamma(\alpha_1\mathbf{u}),\\
		(P,\mathbf{u},S)|_{t=0}=(P_0,\mathbf{u}_0,S_0).
	\end{array}
	\right.
\end{equation}

The existence result in \cite{peng} is in the following.
\begin{theorem}\cite{peng}\label{Theorem 1.0}
	For any $(P_0^l,\mathbf{u}_0^l,S_0^l)\in\dot{B}_{2,1}^{\frac{n}{2}-1}$ and $(\Lambda P_0^h, \mathbf{u}_0^h,\Lambda S_0^h)\in\dot{B}_{2,1}^{\frac{n}{2}+1}$ when $n=2,3$, if there exists a small constant $c_0$ such that
	\begin{equation}\label{1.6}
		\|(P_0^l,\mathbf{u}_0^l,S_0^l)\|_{\dot{B}_{2,1}^{\frac{n}{2}-1}}+\|(\Lambda P_0^h,\mathbf{u}_0^h,\Lambda S_0^h)\|_{\dot{B}_{2,1}^{\frac{n}{2}+1}}\leq c_0,
	\end{equation}
	then the Cauchy problem (\ref{1.5(0)}) has a unique global solution $(P,\mathbf{u},S)$ such that
	\begin{equation*}
		\begin{split}
			&P^l\in C_b(\mathbb{R}^+,\dot{B}_{2,1}^{\frac{n}{2}-1})\cap L^1(\mathbb{R}^+,\dot{B}_{2,1}^{\frac{n}{2}+1}),\  P^h\in C_b(\mathbb{R}^+,\dot{B}_{2,1}^{\frac{n}{2}+2})\cap L^1(\mathbb{R}^+,\dot{B}_{2,1}^{\frac{n}{2}+2}),\\
			&\mathbf{u}^l\in C_b(\mathbb{R}^+,\dot{B}_{2,1}^{\frac{n}{2}-1})\cap L^1(\mathbb{R}^+,\dot{B}_{2,1}^{\frac{n}{2}+1}),\ \mathbf{u}^h\in C_b(\mathbb{R}^+,\dot{B}_{2,1}^{\frac{n}{2}+1})\cap L^1(\mathbb{R}^+,\dot{B}_{2,1}^{\frac{n}{2}+3}),\\
			&S^l\in C_b(\mathbb{R}^+,\dot{B}_{2,1}^{\frac{n}{2}-1}),\ S^h\in C_b(\mathbb{R}^+,\dot{B}_{2,1}^{\frac{n}{2}+2}).
		\end{split}
	\end{equation*}
\end{theorem}

The higher regularity in Theorem \ref{Theorem 1.0} for (\ref{1.5(0)}) than the full Navier-Stokes equation with $\kappa>0$ in \cite{danchin4} is arising from the loss of dissipation of the entropy $s$ and the stronger nonlinear term $\Gamma(\alpha_1 \mathbf{u})\sim|\nabla \mathbf{u}|^2$. In fact, as stated in \cite{peng}, by using classical product law for this nonlinear term, one can easily see that it requires $\mathbf{u}\in\dot{B}_{2,1}^{\frac{n}{2}+1}$ in high frequency, and hence $\mathbf{u}\in L^1(\dot{B}_{2,1}^{\frac{n}{2}+3})$ in the high frequency.

In this article, we aim to construct global solution in $L^p$-framework for $n$-dimensional Cauchy problem (\ref{1.5(0)}) with $n\geq2$, which is a natural generalization of the result in \cite{peng}. Additionally, we established the decay rate of the solution. Due to the wave operator in the low frequency, one can only expect to extend $L^2$-estimate to $L^p$-estimate in the high frequency. Moreover, due to non-dissipative entropy $s$, we can only close the energy estimate for $s$ in $L^2$-norm in the high frequency. The main observation in $L^p$-framework in this article is that we can relax the regularity of $S$ in the high frequency, that is, $S^h\in C_b(\mathbb{R}^+,\dot{B}_{2,1}^{\frac{n}{2}+1})$, which is different from that in \cite{peng}. After that, we can close the desired nonlinear estimates. Some ideas and nonlinear $L^p$-estimates based on Bony decomposition are from those in \cite{xu2023}, where they considered the viscous liquid-gas two-phase flow model.

\noindent{\rm\bf Notations: }
The letter $C$ stands for a generic positive constant whose meaning is clear from the context.
We write $f\lesssim g$ instead of $f\leq Cg$. For operators $A$ and $B$, we denote the commutator $[A,B]=AB-BA$.

Now, the main results on (\ref{1.5(0)}) are stated as follows. In the following, we denote $\dot{B}_{p,r}^\sigma:=\dot{B}_{p,r}^\sigma(\mathbb{R}^n)$ without confusion.
\begin{theorem}[Global existence]\label{Theorem 1.1}
	Let $n\geq2$ and $p$ fulfills
	\begin{equation}\label{1.6}
		2\leq p\leq\min(4,2n/(n-2))\ and, \ additionally,\
		p\neq4\ if\ n=2.
	\end{equation}
	For any $(P_0^l,\mathbf{u}_0^l,S_0^l)\in\dot{B}_{2,1}^{\frac{n}{2}-1}$, $S_0^h\in\dot{B}_{2,1}^{\frac{n}{2}+1}$ and $(\Lambda P_0^h,\mathbf{u}_0^h)\in\dot{B}_{p,1}^{\frac{n}{p}+1}$, if there exists a small constant $c_0$ such that
	\begin{equation}\label{1.7}
		\|(P_0^l,\mathbf{u}_0^l,S_0^l)\|_{\dot{B}_{2,1}^{\frac{n}{2}-1}}+\|S_0^h\|_{\dot{B}_{2,1}^{\frac{n}{2}+1}}+\|(\Lambda P_0^h,\mathbf{u}_0^h)\|_{\dot{B}_{p,1}^{\frac{n}{p}+1}}\leq c_0,
	\end{equation}
	then the Cauchy problem (\ref{1.5(0)}) has a unique global solution $(P,\mathbf{u},S)$ such that
	\begin{equation*}
		\begin{split}
			&P^l\in C_b(\mathbb{R}^+,\dot{B}_{2,1}^{\frac{n}{2}-1})\cap L^1(\mathbb{R}^+,\dot{B}_{2,1}^{\frac{n}{2}+1}),\ \Lambda P^h\in C_b(\mathbb{R}^+,\dot{B}_{p,1}^{\frac{n}{p}+1})\cap L^1(\mathbb{R}^+,\dot{B}_{p,1}^{\frac{n}{p}+1}),\\
			&\mathbf{u}^l\in C_b(\mathbb{R}^+,\dot{B}_{2,1}^{\frac{n}{2}-1})\cap L^1(\mathbb{R}^+,\dot{B}_{2,1}^{\frac{n}{2}+1}),\ \mathbf{u}^h\in C_b(\mathbb{R}^+,\dot{B}_{p,1}^{\frac{n}{p}+1})\cap L^1(\mathbb{R}^+,\dot{B}_{p,1}^{\frac{n}{p}+3}),\\
			&S^l\in C_b(\mathbb{R}^+,\dot{B}_{2,1}^{\frac{n}{2}-1}),\ S^h\in C_b(\mathbb{R}^+,\dot{B}_{2,1}^{\frac{n}{2}+1}).
		\end{split}
	\end{equation*}
	Moreover, there exists some constant $C>0$ such that
	\begin{equation*}\label{1.7}
		\mathcal{X}(t)\leq Cc_0,
	\end{equation*}
	with
	\begin{equation*}
		\begin{split}
			\mathcal{X}(t)\triangleq&\|(P^l,\mathbf{u}^l,S^l)\|_{\tilde{L}_t^\infty(\dot{B}_{2,1}^{\frac{n}{2}-1})}+\alpha_3\|(P^l,\mathbf{u}^l)\|_{L_t^1(\dot{B}_{2,1}^{\frac{n}{2}+1})}\\
			&+\|(\Lambda P^h,\mathbf{u}^h)\|_{\tilde{L}_t^\infty(\dot{B}_{p,1}^{\frac{n}{p}+1})}+\|S^h\|_{\tilde{L}_t^\infty(\dot{B}_{2,1}^{\frac{n}{2}+1})}\\
			&+\|\Lambda P^h\|_{L_t^1(\dot{B}_{p,1}^{\frac{n}{p}+1})}+\|\mathbf{u}^h\|_{L_t^1(\dot{B}_{p,1}^{\frac{n}{p}+3})}.
		\end{split}
	\end{equation*}
\end{theorem}


\begin{remark}\label{Remark 1.1} Compared with the global existence in Theorem 1.1, we not only relax the regularity of $S$ in high frequency, but also extend $L^2$-estimate of $(P,\mathbf{u})$ to $L^p$-estimate of $(P,\mathbf{u})$ with $p$ satisfies (\ref{1.6}) in high frequency. In fact,
	$S^h\in {\tilde{L}_t^\infty(\dot{B}_{2,1}^{\frac{n}{2}+1})}$ is crucial to get the desired estimates for $(\Lambda P^h,\mathbf{u}^h)\|_{\tilde{L}_t^\infty(\dot{B}_{p,1}^{\frac{n}{p}+1})}$, see the details in Sections 2.4-2.6.
\end{remark}

The optimal time-decay estimates of solutions are obtained as follows.

\begin{theorem}\label{Theorem 1.2}
	Let $(P,\mathbf{u},S)$ be the corresponding global solution of (\ref{1.5(0)}) for $n\geq2$. Assume that the real number $\sigma$ satisfies
	\begin{equation}\label{1.10}
		\frac{n}{2}-\frac{2n}{p}\leq\sigma<\frac{n}{2}-1.
	\end{equation}
	If $\|(P_0^l,\mathbf{u}_0^l,S_0^l)\|_{\dot{B}_{2,\infty}^{\sigma}}$ is bounded, then it holds that
	\begin{equation}\label{1.11}
		\|\Lambda^\beta(P,\mathbf{u})\|_{L^p}\lesssim(1+t)^{-\frac{n}{2}(\frac{1}{2}-\frac{1}{p})-\frac{\beta-\sigma}{2}},\ if\ -\tilde{\sigma}_1<\beta\leq\frac{n}{p}-1.
	\end{equation}
	Additionally, when $n\geq3$, one can further get the decay rate for the temporal-derivative of non-dissipative variable $S$ as follows:
	\begin{equation}\label{1.11}
		\|\Lambda^q S_t\|_{L^p}\lesssim(1+t)^{-\frac{n}{2}(\frac{1}{2}-\frac{1}{p})-\frac{q-\sigma+1}{2}},\ if\ -\min\big(\frac{n}{p},\tilde{\sigma}_1+1\big)<q\leq \frac{n}{p}-2,\ p<n,
	\end{equation}
	where $\tilde{\sigma}_1\triangleq-\sigma+n(\frac{1}{2}-\frac{1}{p})$.
\end{theorem}

\begin{remark}The variable $S$ is non-dissipative, but we can get the $L^p$-decay rate of $S_t$ when $n\geq3$, which has not been given in \cite{peng}. Of course, if giving higher regularity of the initial data, one can also get the $L^p$-decay rate for the higher order derivatives of the solution.
\end{remark}

The remainder of the paper is organized as follows.  In Section 2, we mainly give the a $priori$ estimate of the solution in $L^p$-framework, which together with local existence gives the global existence. In Section 3, we derive the $L^p$-decay rate of the solution given in Section 2 under the initial assumption in some negative Besov space. Finally, we recall the Littlewood-Paley decomposition, Besov spaces and related analysis tools in Appendix.

\section{The Proof of Global Existence}

\quad\quad In this section, we will divide into six subsections to verify the global existence in Theorem \ref{Theorem 1.1}. In the first to third subsection, we show the estimate in the low frequency. In the forth to sixth subsection we present the estimate in the high frequency. In the end of this section, we eventually acquire the local existence result in system (\ref{2.1}) and follow the local existence theorem to further prove the global existence in Theorem \ref{Theorem 1.1}.

Recall the system \eqref{1.5(0)}:
\begin{equation}\label{2.1}
	\left\{
	\begin{array}{lc}
		\partial_tP+\alpha_2{\rm div}{\mathbf{u}}=f_1,\\[2mm]
		\partial_t\mathbf{u}-\alpha_3\Delta \mathbf{u}-\alpha_4\nabla{\rm div}\mathbf{u}+\alpha_2\nabla P=f_2,\\[2mm]
		\partial_tS+\alpha_1\mathbf{u}\cdot\nabla S=f_3,\\[2mm]
	\end{array}
\right.
\end{equation}
where
\begin{equation}\label{2.2}
	\begin{split}
		f_1:=&-\alpha_1\mathbf{u}\cdot\nabla P-\frac{C_v+R}{C_v}\alpha_1P{\rm div}\mathbf{u}+\frac{C_v+R}{C_v^2}\Gamma(\alpha_1\mathbf{u}),\\
		f_2:=&-\alpha_1\mathbf{u}\cdot\nabla \mathbf{u}-\frac{1}{\alpha_1}k(\rho)\nabla P+k(\rho)(\mu\Delta \mathbf{u}+(\lambda+\mu)\nabla{\rm div}\mathbf{u}),\\
		f_3:=&\frac{R}{\bar{p}}(1-I(P))\Gamma(\alpha_1\mathbf{u}).
	\end{split}
\end{equation}
Besides, for convenience, we define
\begin{equation*}
	\begin{split}
		\mathscr{E}_\infty(t):=&||(P,\mathbf{u},S)||^l_{\dot{B}^{\frac{n}{2}-1}_{2,1}}+||(\Lambda P, \mathbf{u})||^h_{\dot{B}^{\frac{n}{p}+1}_{p,1}}+||S||^h_{\dot{B}^{\frac{n}{2}+1}_{2,1}},\\
		\mathscr{E}_1(t):=&||(P,\mathbf{u})||^l_{\dot{B}^{\frac{n}{2}+1}_{2,1}}+||\Lambda P||^h_{\dot{B}^{\frac{n}{p}+1}_{p,1}}+||\mathbf{u}||^h_{\dot{B}^{\frac{n}{p}+3}_{p,1}}.		
	\end{split}
\end{equation*}

\subsection{The estimate of $\mathbb{P}\mathbf{u}$ in the low frequency.}
\quad\quad Firstly, we estimate the $\mathbb{P}\mathbf{u}$ in the low frequency. We handle the second equation in \eqref{2.1} by taking the operator  $\mathbb{P}$ and get
\begin{equation*}
	\partial_t\mathbb{P}\mathbf{u}-\alpha_3\Delta\mathbb{P}\mathbf{u}=\mathbb{P}f_2.
\end{equation*}
Applying $\dot{\Delta}_j$ to the above equation, and taking the $L^2$ inner product with $\dot{\Delta}_j\mathbb{P}\mathbf{u}$ gives
\begin{equation}\label{2.3}
	\frac{1}{2}\frac{d}{dt}||\dot{\Delta}_j\mathbb{P}\mathbf{u}||^2_{L^2}+c\alpha_32^{2j}||\dot{\Delta}_j\mathbb{P}\mathbf{u}||^2_{L^2}\lesssim||\dot{\Delta}_j\mathbb{P}f_2||_{L^2}||\dot{\Delta}_j\mathbb{P}\mathbf{u}||^2_{L^2},
\end{equation}
where the term $c2^{2j}||\dot{\Delta}_j\mathbb{P}||^2_{L^2}$ is produced by Bernstein's inequality: there exists a positive constant $c$ so that
\begin{equation*}
	-\int_{\mathbb{R}^n}\Delta\dot{\Delta}_j\mathbb{P}\mathbf{u}\cdot \dot{\Delta}_j\mathbb{P}dx\ge c2^{2j}||\dot{\Delta}_j\mathbb{P}\mathbf{u}||^2_{L^2}.
\end{equation*}
After that, multiplying by $1/||\dot{\Delta}_j\mathbb{P}\mathbf{u}||_{L^2}2^{j(\frac{n}{2}-1)}$ on both hand side of \eqref{2.3} and integrating the inequality from 0 to $t$, we can acquire by summing up about $j\le j_0$ that
\begin{equation}\label{2.4}
	||\mathbb{P}\mathbf{u}^l||_{\widetilde{L}^\infty_{t}(\dot{B}^{\frac{n}{2}-1}_{2,1})}
+\alpha_3||\mathbb{P}\mathbf{u}^l||_{L^1_{t}(\dot{B}^{\frac{n}{2}+1}_{2,1})}
	\lesssim
	||\mathbb{P}\mathbf{u}^l_0||_{\dot{B}^{\frac{n}{2}-1}_{2,1}}+||(f_2)^l||_{L^1_{t}(\dot{B}^{\frac{n}{2}-1}_{2,1})}.
\end{equation}
Now, we start to handle the terms in $f_2$ one by one. From the expression of $f_2$, we have
\begin{equation}\label{2.5}
	\begin{split}
			||{(f_2)}^l||_{\dot{B}^{\frac{n}{2}-1}_{2,1}}\lesssim&||(\mathbf{u}\cdot\nabla \mathbf{u})^l||_{\dot{B}^{\frac{n}{2}-1}_{2,1}}+||(k(\rho)\nabla P)^l||_{\dot{B}^{\frac{n}{2}-1}_{2,1}}\\&+||(k(\rho)\Delta \mathbf{u})^l||_{\dot{B}^{\frac{n}{2}-1}_{2,1}}+||(k(\rho)\nabla{\rm div}\mathbf{u})^l||_{\dot{B}^{\frac{n}{2}-1}_{2,1}}.
	\end{split}
\end{equation}
 In terms of $||(\mathbf{u}\cdot\nabla \mathbf{u})^l||_{\dot{B}^{\frac{n}{2}-1}_{2,1}}$, using Bony decomposition, we have
\begin{equation*}
	\mathbf{u}\cdot \nabla \mathbf{u}=T_\mathbf{u}\nabla u +R(u,\nabla \mathbf{u})+T_{\nabla \mathbf{u}}\mathbf{u}.
\end{equation*}
According to embedding theorem of Proposition \ref{A.5}, interpolation theorem of Proposition \ref{interpolation} and Proposition \ref{A.11}, we can get
\begin{equation}\label{2.6}
	\begin{split}
		||T_{\mathbf{u}}\nabla \mathbf{u}+R(\mathbf{u},\nabla \mathbf{u})||_{\dot{B}^{\frac{n}{2}-1}_{2,1}}
		\lesssim&||\mathbf{u}||_{\dot{B}^{\frac{n}{p}-1}_{p,1}}||\nabla \mathbf{u}||_{\dot{B}^{\frac{n}{p}}_{p,1}}
		\lesssim||\mathbf{u}||_{\dot{B}^{\frac{n}{p}-1}_{p,1}}||\mathbf{u}||_{\dot{B}^{\frac{n}{p}+1}_{p,1}}\\
		\lesssim&(||\mathbf{u}^l||_{\dot{B}^{\frac{n}{p}-1}_{p,1}}
+||\mathbf{u}^h||_{\dot{B}^{\frac{n}{p}+1}_{p,1}})(||\mathbf{u}^l||_{\dot{B}^{\frac{n}{p}+1}_{p,1}}+||\mathbf{u}^h||_{\dot{B}^{\frac{n}{p}+3}_{p,1}})\\
		\lesssim&(||\mathbf{u}^l||_{\dot{B}^{\frac{n}{2}-1}_{2,1}}
+||\mathbf{u}^h||_{\dot{B}^{\frac{n}{p}+1}_{p,1}})(||\mathbf{u}^l||_{\dot{B}^{\frac{n}{2}+1}_{2,1}}+||\mathbf{u}^h||_{\dot{B}^{\frac{n}{p}+3}_{p,1}})\\
		\lesssim&\mathscr{E}_\infty(t)\mathscr{E}_1(t),
	\end{split}
\end{equation}
and
\begin{equation}\label{2.7}
	\begin{split}
		||T_{\nabla \mathbf{u}}\mathbf{u}||_{\dot{B}^{\frac{n}{2}-1}_{2,1}}
		\lesssim&||\nabla \mathbf{u}||_{\dot{B}^{\frac{n}{p}-1}_{p,1}}||\mathbf{u}||_{\dot{B}^{\frac{n}{p}}_{p,1}}
		\lesssim||\mathbf{u}||^2_{\dot{B}^{\frac{n}{p}}_{p,1}}
		\lesssim||\mathbf{u}||_{\dot{B}^{\frac{n}{p}-1}_{p,1}}||\mathbf{u}||_{\dot{B}^{\frac{n}{p}+1}_{p,1}}\\
		\lesssim&(||\mathbf{u}^l||_{\dot{B}^{\frac{n}{2}-1}_{2,1}}
+||\mathbf{u}^h||_{\dot{B}^{\frac{n}{p}+1}_{p,1}})(||\mathbf{u}^l||_{\dot{B}^{\frac{n}{2}+1}_{2,1}}+||\mathbf{u}^h||_{\dot{B}^{\frac{n}{p}+3}_{p,1}})\\
		\lesssim&\mathscr{E}_\infty(t)\mathscr{E}_1(t).
	\end{split}
\end{equation}
Combining with \eqref{2.6} and \eqref{2.7} directly implies that
\begin{equation}\label{2.8}
	||(\mathbf{u}\cdot \nabla \mathbf{u})^l||_{\dot{B}^{\frac{n}{2}-1}_{2,1}}\lesssim\mathscr{E}_\infty(t)\mathscr{E}_1(t).
\end{equation}
In terms of $||(k(\rho)\nabla P)^l||_{\dot{B}^{\frac{n}{2}-1}_{2,1}}$, we also use Bony decomposition and write
\begin{equation*}
	k(\rho)\nabla P=T_{k(\rho)}\nabla P+R(k(\rho),\nabla P)+T_{\nabla P}k(\rho).
\end{equation*}
For the first two terms of decomposition, we need the fact that
\begin{equation*}
	\frac{1}{\rho}-\frac{1}{\bar{\rho}}\sim \mathcal{O}(1)(P+S),
\end{equation*}
so we can acquire
\begin{equation}\label{2.9}
	\begin{split}
		&||T_{k(\rho)}\nabla P+R(k(\rho),\nabla P)||_{\dot{B}^{\frac{n}{2}-1}_{2,1}}
		\lesssim||k(\rho)||_{\dot{B}^{\frac{n}{p}-1}_{p,1}}||\nabla P||_{\dot{B}^{\frac{n}{p}}_{p,1}}\\
		\lesssim&\ ||\rho||_{\dot{B}^{\frac{n}{p}-1}_{p,1}}|| P||_{\dot{B}^{\frac{n}{p}+1}_{p,1}}
		\lesssim||(P,S)||_{\dot{B}^{\frac{n}{p}-1}_{p,1}}|| P||_{\dot{B}^{\frac{n}{p}+1}_{p,1}}\\
		\lesssim&\ (||P^l||_{\dot{B}^{\frac{n}{p}-1}_{p,1}}+||P^h||_{\dot{B}^{\frac{n}{p}+2}_{p,1}}+||S^l||_{\dot{B}^{\frac{n}{p}-1}_{p,1}}
+||S^h||_{\dot{B}^{\frac{n}{p}+1}_{p,1}})\\
		&\ \times(||P^l||_{\dot{B}^{\frac{n}{p}+1}_{p,1}}+||P^h||_{\dot{B}^{\frac{n}{p}+2}_{p,1}})\\
		\lesssim&\ (||P^l||_{\dot{B}^{\frac{n}{2}-1}_{2,1}}+||P^h||_{\dot{B}^{\frac{n}{p}+2}_{p,1}}+||S^l||_{\dot{B}^{\frac{n}{2}-1}_{2,1}}
+||S^h||_{\dot{B}^{\frac{n}{2}+1}_{2,1}})\\
		&\times(||P^l||_{\dot{B}^{\frac{n}{2}+1}_{2,1}}+||P^h||_{\dot{B}^{\frac{n}{p}+2}_{p,1}})\\
		\lesssim&\ \mathscr{E}_\infty(t)\mathscr{E}_1(t).
	\end{split}
\end{equation}
Before dealing with the last terms of  $||T_{\nabla P}k(\rho)||_{\dot{B}^{\frac{n}{2}-1}_{2,1}}$, we first note the fact that
\begin{equation}\label{2.10}
	k(\rho)=k'(0)\rho+\rho\widetilde{k}(\rho)\ \ {\rm with}\ \widetilde{k}(0)=0.
\end{equation}
Hence, we can get
\begin{equation}\label{2.11}
	||T_{\nabla P}k(\rho)||_{\dot{B}^{\frac{n}{2}-1}_{2,1}}\lesssim||T_{\nabla P}\rho||_{\dot{B}^{\frac{n}{2}-1}_{2,1}}+||T_{\nabla P}\rho\widetilde{k}(\rho)||_{\dot{B}^{\frac{n}{2}-1}_{2,1}}.
\end{equation}
Besides, we know $\rho=\rho(P,S)$ and $\rho(0,0)=0$, thus it can be written as
\begin{equation}\label{2.12}
	\rho=\rho'_p(0)P+\rho'_S(0)S+\widetilde{\rho}(P,S)PS\ \ {\rm with}\ \widetilde{\rho}(0,0)=0.
\end{equation}
Using the above equation into the term of $||T_{\nabla P}\rho||_{\dot{B}^{\frac{n}{2}-1}_{2,1}}$ in \eqref{2.11} yields that
\begin{equation}\label{2.13}
	||T_{\nabla P}\rho||_{\dot{B}^{\frac{n}{2}-1}_{2,1}}\lesssim||T_{\nabla P}P||_{\dot{B}^{\frac{n}{2}-1}_{2,1}}+||T_{\nabla P}S||_{\dot{B}^{\frac{n}{2}-1}_{2,1}}+||T_{\nabla P}\widetilde{\rho}(P,S)PS||_{\dot{B}^{\frac{n}{2}-1}_{2,1}}.
\end{equation}
For the first term of the right-hand side of \eqref{2.13}, we have
\begin{equation}\label{2.14}
	\begin{split}
		||T_{\nabla P}P||_{\dot{B}^{\frac{n}{2}-1}_{2,1}}\lesssim
		&\ ||\nabla P||_{\dot{B}^{\frac{n}{p}-1}_{p,1}}||P||_{\dot{B}^{\frac{n}{p}}_{p,1}}
		\lesssim||P||^2_{\dot{B}^{\frac{n}{p}}_{p,1}}
		\lesssim||P||_{\dot{B}^{\frac{n}{p}-1}_{p,1}}||P||_{\dot{B}^{\frac{n}{p}+1}_{p,1}}\\
		\lesssim&\ (||P^l||_{\dot{B}^{\frac{n}{p}-1}_{p,1}}+||P^h||_{\dot{B}^{\frac{n}{p}+2}_{p,1}})(||P^l||_{\dot{B}^{\frac{n}{p}+1}_{p,1}}
+||P^h||_{\dot{B}^{\frac{n}{p}+2}_{p,1}})\\
		\lesssim&\ (||P^l||_{\dot{B}^{\frac{n}{2}-1}_{2,1}}+||P^h||_{\dot{B}^{\frac{n}{p}+2}_{p,1}})(||P^l||_{\dot{B}^{\frac{n}{2}+1}_{2,1}}
+||P^h||_{\dot{B}^{\frac{n}{p}+2}_{p,1}})\\
		\lesssim&\ \mathscr{E}_\infty(t)\mathscr{E}_1(t).
	\end{split}
\end{equation}
For the second term, we have
\begin{equation}\label{2.15}
	\begin{split}
		||T_{\nabla P}S||_{\dot{B}^{\frac{n}{2}-1}_{2,1}}
		\lesssim&\ ||\nabla P^l||_{\dot{B}^{\frac{n}{2}}_{2,1}}||S||_{\dot{B}^{\frac{n}{2}-1}_{2,1}}
		\lesssim||P^l||_{\dot{B}^{\frac{n}{2}+1}_{2,1}}||S||_{\dot{B}^{\frac{n}{2}-1}_{2,1}}\\
		\lesssim&\ ||P^l||_{\dot{B}^{\frac{n}{2}+1}_{2,1}}(||S^l||_{\dot{B}^{\frac{n}{2}-1}_{2,1}}+||S^h||_{\dot{B}^{\frac{n}{2}+1}_{2,1}})\\
		\lesssim&\ \mathscr{E}_1(t)\mathscr{E}_\infty(t).
	\end{split}
\end{equation}
And for the last term, we use Proposition \ref{interpolation} to get
\begin{equation}\label{2.16}
	\begin{split}
		||T_{\nabla P}\widetilde{\rho}(P,S)PS||_{\dot{B}^{\frac{n}{2}-1}_{2,1}}
		\lesssim&\ ||\nabla P||_{\dot{B}^{\frac{n}{p}-1}_{p,1}}||\widetilde{\rho}(P,S)PS||_{\dot{B}^{\frac{n}{p}}_{p,1}}\\
		\lesssim&\ ||P||_{\dot{B}^{\frac{n}{p}}_{p,1}}||(P,S)||_{\dot{B}^{\frac{n}{p}}_{p,1}}||P||_{\dot{B}^{\frac{n}{p}}_{p,1}}||S||_{\dot{B}^{\frac{n}{p}}_{p,1}}\\
		\lesssim&\ ||P||_{\dot{B}^{\frac{n}{p}-1}_{p,1}}||(P,S)||_{\dot{B}^{\frac{n}{p}}_{p,1}}||P||_{\dot{B}^{\frac{n}{p}+1}_{p,1}}||S||_{\dot{B}^{\frac{n}{p}}_{p,1}}\\
		\lesssim&\ \Big(||P^l||_{\dot{B}^{\frac{n}{2}-1}_{2,1}}+||P^h||_{\dot{B}^{\frac{n}{p}+2}_{p,1}}+||S^l||_{\dot{B}^{\frac{n}{2}-1}_{2,1}}+||S^h||_{\dot{B}^{\frac{n}{2}+1}_{2,1}}\Big)^3\\
		&\ \times(||P^l||_{\dot{B}^{\frac{n}{2}+1}_{2,1}}+||P^h||_{\dot{B}^{\frac{n}{p}+2}_{p,1}})\\
		\lesssim&\ (\mathscr{E}_\infty(t))^3\mathscr{E}_1(t).
	\end{split}
\end{equation}
Collecting \eqref{2.14}-\eqref{2.16} into \eqref{2.13}, one has
\begin{equation}\label{2.17}
	||T_{\nabla P}\rho||_{\dot{B}^{\frac{n}{2}-1}_{2,1}}\lesssim(1+(\mathscr{E}_\infty(t))^2)\mathscr{E}_\infty(t)\mathscr{E}_1(t).
\end{equation}
To deal with the term $||T_{\nabla P}\rho\widetilde{k}(\rho)||_{\dot{B}^{\frac{n}{2}-1}_{2,1}}$ by using \eqref{a.15} of Proposition \ref{A.11}, we first have
\begin{equation}\label{2.18}
	||T_{\nabla P}\rho\widetilde{k}(\rho)||_{\dot{B}^{\frac{n}{2}-1}_{2,1}}
	\lesssim||\nabla P^l||_{\dot{B}^{\frac{n}{2}}_{2,1}}||\rho\widetilde{k}(\rho)||_{\dot{B}^{\frac{n}{2}-1}_{2,1}}.
\end{equation}
In terms of $||\rho\widetilde{k}(\rho)||_{\dot{B}^{\frac{n}{2}-1}_{2,1}}$, we use Bony decomposition and Proposition \ref{A.11} again to get
\begin{equation}\label{2.19}
	\begin{split}
		||\rho\widetilde{k}(\rho)||_{\dot{B}^{\frac{n}{2}-1}_{2,1}}\lesssim&\ ||T_\rho\widetilde{k}(\rho)+R(\rho, \widetilde{k}(\rho))||_{\dot{B}^{\frac{n}{2}-1}_{2,1}}+||T_{\widetilde{k}(\rho)}\rho||_{\dot{B}^{\frac{n}{2}-1}_{2,1}}\\
		\lesssim&\ ||\rho||_{\dot{B}^{\frac{n}{p}-1}_{p,1}}||\widetilde{k}(\rho)||_{\dot{B}^{\frac{n}{p}}_{p,1}}+||\widetilde{k}(\rho)||_{\dot{B}^{\frac{n}{p}-1}_{p,1}}||\rho||_{\dot{B}^{\frac{n}{p}}_{p,1}}\\
		\lesssim&\ ||(P,S)||_{\dot{B}^{\frac{n}{p}-1}_{p,1}}||(P,S)||_{\dot{B}^{\frac{n}{p}}_{p,1}}\\
		\lesssim&\ (||P^l||_{\dot{B}^{\frac{n}{2}-1}_{2,1}}+||P^h||_{\dot{B}^{\frac{n}{p}+2}_{p,1}}+||S^l||_{\dot{B}^{\frac{n}{2}-1}_{2,1}}+||S^h||_{\dot{B}^{\frac{n}{2}+1}_{2,1}})^2\\
		\lesssim&\ (\mathscr{E}_\infty(t))^2.
	\end{split}
\end{equation}
So, taking \eqref{2.19} into \eqref{2.18}, one can obtain
\begin{equation}\label{2.20}
	\begin{split}
		||T_{\nabla P}\rho\widetilde{k}(\rho)||_{\dot{B}^{\frac{n}{2}-1}_{2,1}}
		\lesssim&\ ||\nabla P^l||_{\dot{B}^{\frac{n}{2}}_{2,1}}||\rho\widetilde{k}(\rho)||_{\dot{B}^{\frac{n}{2}-1}_{2,1}}
		\lesssim||P^l||_{\dot{B}^{\frac{n}{2}+1}_{2,1}}(\mathscr{E}_\infty(t))^2\\
		\lesssim&\ (\mathscr{E}_\infty(t))^2\mathscr{E}_1(t).
	\end{split}
\end{equation}
Plugging \eqref{2.17} and \eqref{2.20} into \eqref{2.11} gives
\begin{equation}\label{2.21}
	||T_{\nabla P}k(\rho)||_{\dot{B}^{\frac{n}{2}-1}_{2,1}}\lesssim(1+\mathscr{E}_\infty(t)+(\mathscr{E}_\infty(t))^2)\mathscr{E}_\infty(t)\mathscr{E}_1(t).
\end{equation}
Hence, combining \eqref{2.9} and \eqref{2.21}, we can get
\begin{equation}\label{2.22}
	||(k(\rho)\nabla P)^l||_{\dot{B}^{\frac{n}{2}-1}_{2,1}}\lesssim(1+\mathscr{E}_\infty(t)+(\mathscr{E}_\infty(t))^2)\mathscr{E}_\infty(t)\mathscr{E}_1(t).
\end{equation}
As for $||(k(\rho)\Delta \mathbf{u})^l||_{\dot{B}^{\frac{n}{2}-1}_{2,1}}$, we continue to use Bony decomposition, Proposition \ref{A.11} and then acquire
\begin{equation}\label{2.23}
	\begin{split}
		||(k(\rho)\Delta \mathbf{u})^l||_{\dot{B}^{\frac{n}{2}-1}_{2,1}}\lesssim&\ ||T_{k(\rho)}\Delta \mathbf{u}+R(k(\rho),\Delta \mathbf{u})||_{\dot{B}^{\frac{n}{2}-1}_{2,1}}+||T_{\Delta \mathbf{u}}k(\rho)||_{\dot{B}^{\frac{n}{2}-1}_{2,1}}\\
		\lesssim&\ ||k(\rho)||_{\dot{B}^{\frac{n}{p}-1}_{p,1}}||\Delta \mathbf{u}||_{\dot{B}^{\frac{n}{p}}_{p,1}}+||\Delta \mathbf{u}||_{\dot{B}^{\frac{n}{p}-1}_{p,1}}||k(\rho)||_{\dot{B}^{\frac{n}{p}}_{p,1}}\\
		\lesssim&\ ||(P,S)||_{\dot{B}^{\frac{n}{p}-1}_{p,1}}||\mathbf{u}||_{\dot{B}^{\frac{n}{p}+2}_{p,1}}
+||\mathbf{u}||_{\dot{B}^{\frac{n}{p}+1}_{p,1}}||(P,S)||_{\dot{B}^{\frac{n}{p}}_{p,1}}\\
		\lesssim&\ (||P^l||_{\dot{B}^{\frac{n}{2}-1}_{2,1}}+||P^h||_{\dot{B}^{\frac{n}{p}+2}_{p,1}}+||S^l||_{\dot{B}^{\frac{n}{2}-1}_{2,1}}+||S^h||_{\dot{B}^{\frac{n}{2}+1}_{2,1}})\\
		&\ \times(||\mathbf{u}^l||_{\dot{B}^{\frac{n}{2}+1}_{2,1}}+||\mathbf{u}^h||_{\dot{B}^{\frac{n}{p}+3}_{p,1}})\\
		\lesssim&\ \mathscr{E}_\infty(t)\mathscr{E}_1(t).
	\end{split}
\end{equation}
Bounding the last terms $||k(\rho)\nabla{\rm div}\mathbf{u}||_{\dot{B}^{\frac{n}{2}-1}_{2,1}}$ is similar, and it also holds that
\begin{equation}\label{2.24}
	||(k(\rho)\nabla{\rm div}\mathbf{u})^l||_{\dot{B}^{\frac{n}{2}-1}_{2,1}}\lesssim\mathscr{E}_\infty(t)\mathscr{E}_1(t).
\end{equation}
So we omit the details here.
Plugging \eqref{2.8}, \eqref{2.22}-\eqref{2.24} into \eqref{2.5}, then $||{(f_2)}^l||_{\dot{B}^{\frac{n}{2}-1}_{2,1}}$ can be estimated as
\begin{equation}\label{2.25}
	||{(f_2)}^l||_{\dot{B}^{\frac{n}{2}-1}_{2,1}}\lesssim(1+\mathscr{E}_\infty(t)+(\mathscr{E}_\infty(t))^2)\mathscr{E}_\infty(t)\mathscr{E}_1(t).
\end{equation}
Consequently, inserting \eqref{2.25} into \eqref{2.4}, we acquire the low frequency $\mathbb{P}u^l$ as
\begin{equation}\label{2.26}
\begin{array}{rl}	&||\mathbb{P}\mathbf{u}^l||_{\widetilde{L}^\infty_{t}(\dot{B}^{\frac{n}{2}-1}_{2,1})}+\alpha_3||\mathbb{P}\mathbf{u}^l||_{L^1_{t}(\dot{B}^{\frac{n}{2}+1}_{2,1})}\\
	\lesssim& \displaystyle ||\mathbb{P}\mathbf{u}^l_0||_{\dot{B}^{\frac{n}{2}-1}_{2,1}}+\int_{0}^{t}(1+\mathscr{E}_\infty(\tau)+(\mathscr{E}_\infty(\tau))^2)\mathscr{E}_\infty(\tau)\mathscr{E}_1(\tau)d\tau.
\end{array}
\end{equation}

\subsection{The estimate of $S$ in the low frequency.}
\quad\quad In this subsection, we will estimate the $S$ in the low frequency. Because there is no dissipation in the equation of $S$, we have introduced the commutator's argument to handle the term $u\cdot\nabla S$. Now, applying $\dot{\Delta}_j$ to the third equation of \eqref{2.1}, we have
\begin{equation*}
	\partial_t\dot{\Delta}_jS+\alpha_1\mathbf{u}\cdot\nabla\dot{\Delta}_jS+[\dot{\Delta}_j,\mathbf{u}\cdot\nabla]S=\dot{\Delta}_jf_3.
\end{equation*}
In the same way as in \eqref{2.4}, we can acquire the derivation by standard energy argument that
\begin{equation}\label{2.27}
	\begin{split}
		||S^l||_{\widetilde{L}^\infty_t(\dot{B}^{\frac{n}{2}-1}_{2,1})}\lesssim&\ ||S_0^l||_{\dot{B}^{\frac{n}{2}-1}_{2,1}}
+||(f_3)^l||_{L_t^1(\dot{B}^{\frac{n}{2}-1}_{2,1})}+\int_{0}^{t}||{\rm div}\mathbf{u}||_{L^\infty}||S^l||_{\dot{B}^{\frac{n}{2}-1}_{2,1}}d\tau\\
		&\ +\int_{0}^{t}\sum\limits_{j\le j_0}2^{(\frac{n}{2}-1)j}||[\dot{\Delta}_j, \mathbf{u}\cdot\nabla]S||_{L^2}d\tau.
	\end{split}
\end{equation}
Following the expression of $f_3$, we have
\begin{equation}\label{2.28}
	||(f_3)^l||_{\dot{B}^{\frac{n}{2}-1}_{2,1}}
\lesssim||(\Gamma(\alpha_1\mathbf{u}))^l||_{\dot{B}^{\frac{n}{2}-1}_{2,1}}+||(I(P)\Gamma(\alpha_1\mathbf{u}))^l||_{\dot{B}^{\frac{n}{2}-1}_{2,1}}.
\end{equation}
For the first term $||(\Gamma(\alpha_1\mathbf{u}))^l||_{\dot{B}^{\frac{n}{2}-1}_{2,1}}$, using Bony decomposition and Proposition \ref{A.11}, it can be bounded by
\begin{equation}\label{2.29}
	\begin{split}
		&\ ||(\Gamma(\alpha_1\mathbf{u}))^l||_{\dot{B}^{\frac{n}{2}-1}_{2,1}}\\
\lesssim&\ ||(\nabla \mathbf{u})^2||_{\dot{B}^{\frac{n}{2}-1}_{2,1}}
		\lesssim||T_{\nabla \mathbf{u}}\nabla \mathbf{u}+R(\nabla \mathbf{u},\nabla \mathbf{u})||_{\dot{B}^{\frac{n}{2}-1}_{2,1}}\\
		\lesssim&\ ||\nabla \mathbf{u}||_{\dot{B}^{\frac{n}{p}-1}_{p,1}}||\nabla \mathbf{u}||_{\dot{B}^{\frac{n}{p}}_{p,1}}
		\lesssim||\mathbf{u}||_{\dot{B}^{\frac{n}{p}}_{p,1}}||\mathbf{u}||_{\dot{B}^{\frac{n}{p}+1}_{p,1}}\\ \lesssim&\ (||\mathbf{u}^l||_{\dot{B}^{\frac{n}{2}-1}_{2,1}}+||\mathbf{u}^h||_{\dot{B}^{\frac{n}{p}+1}_{p,1}})(||\mathbf{u}^l||_{\dot{B}^{\frac{n}{2}+1}_{2,1}}+||\mathbf{u}^h||_{\dot{B}^{\frac{n}{p}+3}_{p,1}})\\
		\lesssim&\ \mathscr{E}_\infty(t)\mathscr{E}_1(t).
	\end{split}
\end{equation}
Next, for the second term $||(I(P)\Gamma(\alpha_1\mathbf{u}))^l||_{\dot{B}^{\frac{n}{2}-1}_{2,1}}$, we use the result of \eqref{2.29}, Proposition \ref{A.7} and Proposition \ref{A.13} to get
\begin{equation}\label{2.30}
	\begin{split}		||(I(P)\Gamma(\alpha_1\mathbf{u}))^l||_{\dot{B}^{\frac{n}{2}-1}_{2,1}}\lesssim&\ ||I(P)||_{\dot{B}^{\frac{n}{p}}_{p,1}}||\Gamma(\alpha_1\mathbf{u})||_{\dot{B}^{\frac{n}{2}-1}_{2,1}}\\
		\lesssim&\ ||P||_{\dot{B}^{\frac{n}{p}}_{p,1}}\mathscr{E}_\infty(t)\mathscr{E}_1(t)\\
		\lesssim&\ (||P^l||_{\dot{B}^{\frac{n}{2}-1}_{2,1}}+||P^h||_{\dot{B}^{\frac{n}{p}+2}_{p,1}})\mathscr{E}_\infty(t)\mathscr{E}_1(t)\\
		\lesssim&\ (\mathscr{E}_\infty(t))^2\mathscr{E}_1(t).
	\end{split}
\end{equation}
Bring \eqref{2.29} and \eqref{2.30} into \eqref{2.28}, one has
\begin{equation}\label{2.31}
	||(f_3)^l||_{\dot{B}^{\frac{n}{2}-1}_{2,1}}\lesssim(1+\mathscr{E}_\infty(t))\mathscr{E}_\infty(t)\mathscr{E}_1(t).
\end{equation}
Thanks to the embedding relation $\dot{B}^{\frac{n}{p}}_{p,1}\hookrightarrow L^\infty$ and \eqref{a.17} of Lemma \ref{A.12}, the last two terms in \eqref{2.27} can be bounded as
\begin{equation}\label{2.32}
	\begin{split}
		&||{\rm div}\mathbf{u}||_{L^\infty}||S^l||_{\dot{B}^{\frac{n}{2}-1}_{2,1}}+\sum\limits_{j\le j_0}2^{(\frac{n}{2}-1)j}||[\dot{\Delta}_j, \mathbf{u}\cdot\nabla]S||_{L^2}\\
		\lesssim\ &||{\rm div}\mathbf{u}||_{\dot{B}^{\frac{n}{p}}_{p,1}}||S^l||_{\dot{B}^{\frac{n}{2}-1}_{2,1}}+||\nabla \mathbf{u}||_{\dot{B}^{\frac{n}{p}}_{p,1}}||S||_{\dot{B}^{\frac{n}{2}-1}_{2,1}}\\
		\lesssim\ &||\mathbf{u}||_{\dot{B}^{\frac{n}{p}+1}_{p,1}}||S^l||_{\dot{B}^{\frac{n}{2}-1}_{2,1}}+||\mathbf{u}||_{\dot{B}^{\frac{n}{p}+1}_{p,1}}||S||_{\dot{B}^{\frac{n}{2}-1}_{2,1}}\\
		\lesssim\ &(||\mathbf{u}^l||_{\dot{B}^{\frac{n}{2}+1}_{2,1}}+||\mathbf{u}^h||_{\dot{B}^{\frac{n}{p}+3}_{p,1}})(||S^l||_{\dot{B}^{\frac{n}{2}-1}_{2,1}}+||S^h||_{\dot{B}^{\frac{n}{2}+1}_{2,1}})\\
		\lesssim\ &\mathscr{E}_\infty(t)\mathscr{E}_1(t).
	\end{split}
\end{equation}
As a result, plugging \eqref{2.31}, \eqref{2.32} into \eqref{2.27}, we have
\begin{equation}\label{2.33}
	||S^l||_{\widetilde{L}^\infty_t(\dot{B}^{\frac{n}{2}-1}_{2,1})}\lesssim||S_0^l||_{\dot{B}^{\frac{n}{2}-1}_{2,1}}+\int_{0}^{t}(1+\mathscr{E}_\infty(\tau)+(\mathscr{E}_\infty(\tau))^2)\mathscr{E}_\infty(\tau)\mathscr{E}_1(\tau)d\tau.
\end{equation}

\subsection{The estimate of $(P,\mathbb{Q}\mathbf{u})$ in the low frequency.}
\quad\quad In this subsection, we consider the system of $(P,\mathbb{Q}\mathbf{u})$. Under the barotropic linearized equations (see \cite{Bahouri-2011} or \cite{danchin1}) of energy estimates, we can obtain
\begin{equation}\label{2.34}
	\begin{split}
		&||(P,\mathbb{Q}\mathbf{u})^l||_{\widetilde{L}^\infty_t(\dot{B}^{\frac{n}{2}-1}_{2,1})}
+(\alpha_3+\alpha_4)||(P,\mathbb{Q}\mathbf{u})^l||_{L^1_{t}(\dot{B}^{\frac{n}{2}+1}_{2,1})}\\
		\lesssim&||(P_0,\mathbb{Q}\mathbf{u}_0)^l||_{\dot{B}^{\frac{n}{2}-1}_{2,1}}+\int_{0}^{t}||(f_1)^l||_{\dot{B}^{\frac{n}{2}-1}_{2,1}}d\tau
+\int_{0}^{t}||(f_2)^l||_{\dot{B}^{\frac{n}{2}-1}_{2,1}}d\tau.
	\end{split}
\end{equation}
It follows from the expression of $||(f_1)^l||_{\dot{B}^{\frac{n}{2}-1}_{2,1}}$ that
\begin{equation}\label{2.35}
	||(f_1)^l||_{\dot{B}^{\frac{n}{2}-1}_{2,1}}\lesssim||(\mathbf{u}\cdot\nabla P)^l||_{\dot{B}^{\frac{n}{2}-1}_{2,1}}+||(P{\rm div}\mathbf{u})^l||_{\dot{B}^{\frac{n}{2}-1}_{2,1}}+||(\Gamma(\alpha_1\mathbf{u}))^l||_{\dot{B}^{\frac{n}{2}-1}_{2,1}}.
\end{equation}
Bony decomposition gives
\begin{equation}\label{2.36}
	\mathbf{u}\cdot\nabla P=T_{\mathbf{u}}\nabla P+R(\mathbf{u},\nabla P)+T_{\nabla P}\mathbf{u}.
\end{equation}
Thus, by using Proposition \ref{A.11}, it can be shown that
\begin{equation}\label{2.37}
	\begin{split}
		||T_{\mathbf{u}}\nabla P+R(\mathbf{u},\nabla P)||_{\dot{B}^{\frac{n}{2}-1}_{2,1}}\lesssim&||\mathbf{u}||_{\dot{B}^{\frac{n}{p}-1}_{p,1}}||\nabla P||_{\dot{B}^{\frac{n}{p}}_{p,1}}
		\lesssim||\mathbf{u}||_{\dot{B}^{\frac{n}{p}-1}_{p,1}}||P||_{\dot{B}^{\frac{n}{p}+1}_{p,1}}\\
		\lesssim&(||\mathbf{u}^l||_{\dot{B}^{\frac{n}{2}-1}_{2,1}}+||\mathbf{u}^h||_{\dot{B}^{\frac{n}{p}+1}_{p,1}})(||P^l||_{\dot{B}^{\frac{n}{2}+1}_{2,1}}+||P^h||_{\dot{B}^{\frac{n}{p}+2}_{p,1}})\\
		\lesssim&\mathscr{E}_\infty(t)\mathscr{E}_1(t),
	\end{split}
\end{equation}
and
\begin{equation}\label{2.38}
	\begin{split}
		||T_{\nabla P}\mathbf{u}||_{\dot{B}^{\frac{n}{2}-1}_{2,1}}\lesssim&||\nabla P||_{\dot{B}^{\frac{n}{p}-1}_{p,1}}||\mathbf{u}||_{\dot{B}^{\frac{n}{p}}_{p,1}}
		\lesssim||P||_{\dot{B}^{\frac{n}{p}}_{p,1}}||\mathbf{u}||_{\dot{B}^{\frac{n}{p}}_{p,1}}\\
		\lesssim&||P||^{\frac{1}{2}}_{\dot{B}^{\frac{n}{p}-1}_{p,1}}||P||^{\frac{1}{2}}_{\dot{B}^{\frac{n}{p}+1}_{p,1}}|
|\mathbf{u}||^{\frac{1}{2}}_{\dot{B}^{\frac{n}{p}-1}_{p,1}}||\mathbf{u}||^{\frac{1}{2}}_{\dot{B}^{\frac{n}{p}+1}_{p,1}}\\
		\lesssim&(||P^l||_{\dot{B}^{\frac{n}{2}-1}_{2,1}}+||P^h||_{\dot{B}^{\frac{n}{p}+2}_{p,1}})^{\frac{1}{2}}(||P^l||_{\dot{B}^{\frac{n}{2}+1}_{2,1}}
+||P^h||_{\dot{B}^{\frac{n}{p}+2}_{p,1}})^{\frac{1}{2}}\\
		&\times(||\mathbf{u}^l||_{\dot{B}^{\frac{n}{2}-1}_{2,1}}
+||\mathbf{u}^h||_{\dot{B}^{\frac{n}{p}+1}_{p,1}})^{\frac{1}{2}}(||\mathbf{u}^l||_{\dot{B}^{\frac{n}{2}+1}_{2,1}}+||\mathbf{u}^h||_{\dot{B}^{\frac{n}{p}+3}_{p,1}})^{\frac{1}{2}}\\
		\lesssim&\mathscr{E}_\infty(t)\mathscr{E}_1(t).
	\end{split}
\end{equation}
From the above inequalities, we can acquire
\begin{equation}\label{2.39}
	||\mathbf{u}\cdot\nabla P||_{\dot{B}^{\frac{n}{2}-1}_{2,1}}\lesssim\mathscr{E}_\infty(t)\mathscr{E}_1(t).
\end{equation}
Similarly, we can use the same way to treat $||P{\rm div}\mathbf{u}||_{\dot{B}^{\frac{n}{2}-1}_{2,1}}$ by using Bony decomposition as
\begin{equation}\label{2.40}
	P{\rm div}\mathbf{u}=T_{P}{\rm div}\mathbf{u}+R(P,{\rm div}\mathbf{u})+T_{{\rm div}\mathbf{u}}P.
\end{equation}
Then Proposition \ref{A.11} leads to
\begin{equation}\label{2.41}
	\begin{split}
		||T_{P}{\rm div}\mathbf{u}+R(P,{\rm div}\mathbf{u})||_{\dot{B}^{\frac{n}{2}-1}_{2,1}}\lesssim&||P||_{\dot{B}^{\frac{n}{p}-1}_{p,1}}||{\rm div}\mathbf{u}||_{\dot{B}^{\frac{n}{p}}_{p,1}}
		\lesssim||P||_{\dot{B}^{\frac{n}{p}-1}_{p,1}}||\mathbf{u}||_{\dot{B}^{\frac{n}{p}+1}_{p,1}}\\
		\lesssim&(||P^l||_{\dot{B}^{\frac{n}{2}-1}_{2,1}}+||P^h||_{\dot{B}^{\frac{n}{p}+2}_{p,1}})(||\mathbf{u}^l||_{\dot{B}^{\frac{n}{2}+1}_{2,1}}+||\mathbf{u}^h||_{\dot{B}^{\frac{n}{p}+3}_{p,1}})\\
		\lesssim&\mathscr{E}_\infty(t)\mathscr{E}_1(t),
	\end{split}
\end{equation}
and
\begin{equation}\label{2.42}
	\begin{split}
		||T_{{\rm div}u}P||_{\dot{B}^{\frac{n}{2}-1}_{2,1}}\lesssim&||{\rm div}\mathbf{u}||_{\dot{B}^{\frac{n}{p}-1}_{p,1}}||P||_{\dot{B}^{\frac{n}{p}}_{p,1}}
		\lesssim||\mathbf{u}||_{\dot{B}^{\frac{n}{p}}_{p,1}}||P||_{\dot{B}^{\frac{n}{p}}_{p,1}}\\
		\lesssim&||\mathbf{u}||^{\frac{1}{2}}_{\dot{B}^{\frac{n}{p}-1}_{p,1}}|
|\mathbf{u}||^{\frac{1}{2}}_{\dot{B}^{\frac{n}{p}+1}_{p,1}}||P||^{\frac{1}{2}}_{\dot{B}^{\frac{n}{p}-1}_{p,1}}||P||^{\frac{1}{2}}_{\dot{B}^{\frac{n}{p}+1}_{p,1}}\\
		\lesssim&(||\mathbf{u}^l||_{\dot{B}^{\frac{n}{2}-1}_{2,1}}+||\mathbf{u}^h||_{\dot{B}^{\frac{n}{p}+1}_{p,1}})^{\frac{1}{2}}
(||\mathbf{u}^l||_{\dot{B}^{\frac{n}{2}+1}_{2,1}}+||\mathbf{u}^h||_{\dot{B}^{\frac{n}{p}+3}_{p,1}})^{\frac{1}{2}}\\
		&\times(||P^l||_{\dot{B}^{\frac{n}{2}-1}_{2,1}}+||P^h||_{\dot{B}^{\frac{n}{p}+2}_{p,1}})^{\frac{1}{2}}(||P^l||_{\dot{B}^{\frac{n}{2}+1}_{2,1}}
+||P^h||_{\dot{B}^{\frac{n}{p}+2}_{p,1}})^{\frac{1}{2}}\\
		\lesssim&\mathscr{E}_\infty(t)\mathscr{E}_1(t),
	\end{split}
\end{equation}
and hence
\begin{equation}\label{2.43}
	||P{\rm div}\mathbf{u}||_{\dot{B}^{\frac{n}{2}-1}_{2,1}}\lesssim\mathscr{E}_\infty(t)\mathscr{E}_1(t).
\end{equation}
The last term $||(\Gamma(\alpha_1\mathbf{u}))^l||_{\dot{B}^{\frac{n}{2}-1}_{2,1}}$ has been verified in \eqref{2.29}, which combining with \eqref{2.39} and \eqref{2.43} into \eqref{2.35} implies
\begin{equation}\label{2.44}
	||(f_1)^l||_{\dot{B}^{\frac{n}{2}-1}_{2,1}}\lesssim\mathscr{E}_\infty(t)\mathscr{E}_1(t).
\end{equation}
Since the term $||(f_2)^l||_{\dot{B}^{\frac{n}{2}-1}_{2,1}}$ has been estimated in \eqref{2.25}, it follows that
\begin{equation}\label{2.45}
	\begin{split}
		&||(P,\mathbb{Q}\mathbf{u})^l||_{\widetilde{L}^\infty_t(\dot{B}^{\frac{n}{2}-1}_{2,1})}+(\alpha_3+\alpha_4)||(P,\mathbb{Q}\mathbf{u})^l||_{L^1_{t}(\dot{B}^{\frac{n}{2}+1}_{2,1})}\\
		\lesssim&||(P_0,\mathbb{Q}\mathbf{u}_0^l)||_{\dot{B}^{\frac{n}{2}-1}_{2,1}}+\int_{0}^{t}(1+\mathscr{E}_\infty(\tau)+(\mathscr{E}_\infty(\tau))^2)\mathscr{E}_\infty(\tau)\mathscr{E}_1(\tau)d\tau.
	\end{split}
\end{equation}

\subsection{The estimate of $\mathbb{P}\mathbf{u}$ in the high frequency.}
\quad\quad Next, we are ready to estimate the terms of $\mathbb{P}\mathbf{u}$ in the high frequency. The same as \eqref{2.4} in Subsection 2.1, we take the operator $\mathbb{P}$ again to the second equation in \eqref{2.1} and get
\begin{equation*}
	\partial_t\mathbb{P}\mathbf{u}-\alpha_3\Delta\mathbb{P}\mathbf{u}=\mathbb{P}f_2.
\end{equation*}
Applying $\dot{\Delta}_j$ to the above equation, after that taking the $L^2$ inner product with $|\dot{\Delta}_j\mathbb{P}\mathbf{u}|^{p-2}\dot{\Delta}_j\mathbb{P}\mathbf{u}$ and following the H$\ddot{o}$lder inequality results in
\begin{equation}\label{2.46}
	\frac{1}{p}\frac{d}{dt}||\dot{\Delta}_j\mathbb{P}\mathbf{u}||^p_{L^p}+c\alpha_32^{2j}||\dot{\Delta}_j\mathbb{P}\mathbf{u}||^p_{L^p}
\lesssim||\dot{\Delta}_j\mathbb{P}f_2||_{L^p}||\dot{\Delta}_j\mathbb{P}\mathbf{u}||^{p-1}_{L^p},
\end{equation}
	in which we get the term $c2^{2j}||\dot{\Delta}_j\mathbb{P}\mathbf{u}||^p_{L^p}$ by using the Proposition \ref{A.6}: there exists a positive constant $c$ so that
\begin{equation*}
	-\int_{\mathbb{R}^n}\Delta\dot{\Delta}_j\mathbb{P}\mathbf{u}\cdot|\dot{\Delta}_j\mathbb{P}\mathbf{u}|^{p-2}\dot{\Delta}_j\mathbb{P}\mathbf{u}dx
=(p-1)\int_{\mathbb{R}^n}|\nabla\dot{\Delta}_j\mathbb{P}\mathbf{u}|^2|\dot{\Delta}_j\mathbb{P}\mathbf{u}|^{p-2}dx\ge c2^{2j}||\dot{\Delta}_j\mathbb{P}\mathbf{u}||^p_{L^p}.
\end{equation*}
Afterwards, multiplying by $1/||\dot{\Delta}_j\mathbb{P}\mathbf{u}||^{p-1}_{L^p}2^{j(\frac{n}{p}+1)}$ on both hand side of \eqref{2.46}, integrating the resultant inequality from 0 to $t$, we can obtain by summing up about $j\ge j_0$ that
\begin{equation}\label{2.47}
	||\mathbb{P}\mathbf{u}^h||_{\widetilde{L}^\infty_{t}(\dot{B}^{\frac{n}{p}+1}_{p,1})}
+\alpha_3||\mathbb{P}\mathbf{u}^h||_{L^1_{t}(\dot{B}^{\frac{n}{p}+3}_{p,1})}
	\lesssim
	||\mathbb{P}\mathbf{u}^h_0||_{\dot{B}^{\frac{n}{p}+1}_{p,1}}+||{(f_2)}^h||_{L^1_{t}(\dot{B}^{\frac{n}{p}+1}_{p,1})}.
\end{equation}
In terms of $||{(f_2)}^h||_{\dot{B}^{\frac{n}{p}+1}_{p,1}}$, we have
\begin{equation}\label{2.48}
	\begin{split}
		||{(f_2)}^h||_{\dot{B}^{\frac{n}{p}+1}_{p,1}}\lesssim&||(\mathbf{u}\cdot\nabla \mathbf{u})^h||_{\dot{B}^{\frac{n}{p}+1}_{p,1}}+||(k(\rho)\nabla P)^h||_{\dot{B}^{\frac{n}{p}+1}_{p,1}}\\&+||(k(\rho)\Delta \mathbf{u})^h||_{\dot{B}^{\frac{n}{P}+1}_{p,1}}+||(k(\rho)\nabla{\rm div}\mathbf{u})^h||_{\dot{B}^{\frac{n}{p}+1}_{p,1}}.
	\end{split}
\end{equation}
Now, we start to estimate the right-hand terms in \eqref{2.48} by using \eqref{a.7} of Proposition \ref{A.7}. For the first term $||(\mathbf{u}\cdot\nabla \mathbf{u})^h||_{\dot{B}^{\frac{n}{p}+1}_{p,1}}$, we can get
\begin{equation}\label{2.49}
	\begin{split}
		||(\mathbf{u}\cdot\nabla \mathbf{u})^h||_{\dot{B}^{\frac{n}{p}+1}_{p,1}}\lesssim&||\mathbf{u}||_{L^\infty}||\nabla \mathbf{u}||_{\dot{B}^{\frac{n}{p}+1}_{p,1}}+||\nabla \mathbf{u}||_{L^\infty}||\mathbf{u}||_{\dot{B}^{\frac{n}{p}+1}_{p,1}}\\
		\lesssim&||\mathbf{u}||_{\dot{B}^{\frac{n}{p}}_{p,1}}||\nabla \mathbf{u}||_{\dot{B}^{\frac{n}{p}+1}_{p,1}}+||\nabla \mathbf{u}||_{\dot{B}^{\frac{n}{p}}_{p,1}}||\mathbf{u}||_{\dot{B}^{\frac{n}{p}+1}_{p,1}}\\
		\lesssim&||\mathbf{u}||_{\dot{B}^{\frac{n}{p}}_{p,1}}||\mathbf{u}||_{\dot{B}^{\frac{n}{p}+2}_{p,1}}
+||\mathbf{u}||_{\dot{B}^{\frac{n}{p}+1}_{p,1}}||\mathbf{u}||_{\dot{B}^{\frac{n}{p}+1}_{p,1}}\\
		\lesssim&(||\mathbf{u}^l||_{\dot{B}^{\frac{n}{2}-1}_{2,1}}
+||\mathbf{u}^h||_{\dot{B}^{\frac{n}{p}+1}_{p,1}})(||\mathbf{u}^l||_{\dot{B}^{\frac{n}{2}+1}_{2,1}}+||\mathbf{u}^h||_{\dot{B}^{\frac{n}{p}+3}_{p,1}})\\
		\lesssim&\mathscr{E}_\infty(t)\mathscr{E}_1(t).
	\end{split}
\end{equation}
For the second term $||(k(\rho)\nabla P)^h||_{\dot{B}^{\frac{n}{p}+1}_{p,1}}$, we have
\begin{equation}\label{2.50}
	\begin{split}
		||(k(\rho)\nabla P)^h||_{\dot{B}^{\frac{n}{p}+1}_{p,1}}\lesssim&||k(\rho)||_{L^\infty}||\nabla P||_{\dot{B}^{\frac{n}{p}+1}_{p,1}}+||\nabla P||_{L^\infty}||k(\rho)||_{\dot{B}^{\frac{n}{p}+1}_{p,1}}\\
		\lesssim&||k(\rho)||_{\dot{B}^{\frac{n}{p}}_{p,1}}||\nabla P||_{\dot{B}^{\frac{n}{p}+1}_{p,1}}+||\nabla P||_{\dot{B}^{\frac{n}{p}}_{p,1}}||k(\rho)||_{\dot{B}^{\frac{n}{p}+1}_{p,1}}\\
		\lesssim&||(P,S)||_{\dot{B}^{\frac{n}{p}}_{p,1}}||P||_{\dot{B}^{\frac{n}{p}+2}_{p,1}}+||P||_{\dot{B}^{\frac{n}{p}+1}_{p,1}}||(P,S)||_{\dot{B}^{\frac{n}{p}+1}_{p,1}}\\
		\lesssim&(||P^l||_{\dot{B}^{\frac{n}{2}-1}_{2,1}}+||P^h||_{\dot{B}^{\frac{n}{p}+2}_{p,1}}+||S^l||_{\dot{B}^{\frac{n}{2}-1}_{2,1}}+||S^h||_{\dot{B}^{\frac{n}{2}+1}_{2,1}})\\
		&\times(||P^l||_{\dot{B}^{\frac{n}{2}+1}_{2,1}}+||P^h||_{\dot{B}^{\frac{n}{p}+2}_{p,1}})\\
		\lesssim&\mathscr{E}_\infty(t)\mathscr{E}_1(t).
	\end{split}
\end{equation}
The third and last terms $||(k(\rho)\Delta \mathbf{u})^h||_{\dot{B}^{\frac{n}{P}+1}_{p,1}}$ and $||(k(\rho)\nabla{\rm div}\mathbf{u})^h||_{\dot{B}^{\frac{n}{p}+1}_{p,1}}$ can be bounded as
\begin{equation}\label{2.51}
	\begin{split}
		&||(k(\rho)\Delta \mathbf{u})^h||_{\dot{B}^{\frac{n}{P}+1}_{p,1}}+||(k(\rho)\nabla{\rm div}\mathbf{u})^h||_{\dot{B}^{\frac{n}{p}+1}_{p,1}}
		\lesssim||(k(\rho)\nabla^2\mathbf{u})^h||_{\dot{B}^{\frac{n}{p}+1}_{p,1}}\\
		\lesssim&||k(\rho)||_{L^\infty}||\nabla^2\mathbf{u}||_{\dot{B}^{\frac{n}{p}+1}_{p,1}}
+||\nabla^2\mathbf{u}||_{L^\infty}||k(\rho)||_{\dot{B}^{\frac{n}{p}+1}_{p,1}}\\
		\lesssim&||k(\rho)||_{\dot{B}^{\frac{n}{p}}_{p,1}}||\nabla^2\mathbf{u}||_{\dot{B}^{\frac{n}{p}+1}_{p,1}}
+||\nabla^2\mathbf{u}||_{\dot{B}^{\frac{n}{p}}_{p,1}}||k(\rho)||_{\dot{B}^{\frac{n}{p}+1}_{p,1}}\\
		\lesssim&||(P,S)||_{\dot{B}^{\frac{n}{p}}_{p,1}}||\mathbf{u}||_{\dot{B}^{\frac{n}{p}+3}_{p,1}}
+||\mathbf{u}||_{\dot{B}^{\frac{n}{p}+2}_{p,1}}||(P,S)||_{\dot{B}^{\frac{n}{p}+1}_{p,1}}\\
		\lesssim&(||P^l||_{\dot{B}^{\frac{n}{2}-1}_{2,1}}+||P^h||_{\dot{B}^{\frac{n}{p}+2}_{p,1}}+||S^l||_{\dot{B}^{\frac{n}{2}-1}_{2,1}}
+||S^h||_{\dot{B}^{\frac{n}{2}+1}_{2,1}})
		\times(||\mathbf{u}^l||_{\dot{B}^{\frac{n}{2}+1}_{2,1}}+||\mathbf{u}^h||_{\dot{B}^{\frac{n}{p}+3}_{p,1}})\\
		\lesssim&\mathscr{E}_\infty(t)\mathscr{E}_1(t).
	\end{split}
\end{equation}
So, plugging \eqref{2.49}-\eqref{2.51} into \eqref{2.48}, we get
\begin{equation}\label{2.52}
	||{(f_2)}^h||_{\dot{B}^{\frac{n}{p}+1}_{p,1}}\lesssim\mathscr{E}_\infty(t)\mathscr{E}_1(t).
\end{equation}
Ultimately, it holds that
\begin{equation}\label{2.53}
\begin{split} &||\mathbb{P}\mathbf{u}^h||_{\widetilde{L}^\infty_{t}(\dot{B}^{\frac{n}{p}+1}_{p,1})}+\alpha_3||\mathbb{P}\mathbf{u}^h||_{L^1_{t}(\dot{B}^{\frac{n}{p}+3}_{p,1})}\\
	\lesssim& ||\mathbb{P}\mathbf{u}^h_0||_{\dot{B}^{\frac{n}{p}+1}_{p,1}}
+\int_{0}^{t}(1+\mathscr{E}_\infty(\tau)+(\mathscr{E}_\infty(\tau))^2)\mathscr{E}_\infty(\tau)\mathscr{E}_1(\tau)d\tau.
\end{split}
\end{equation}

\subsection{The estimate of $S$ in the high frequency.}
\quad\quad Now, we turn to estimate the term $S$ in the high frequency. Similar to \eqref{2.27}, we also have brought in the commutator's argument to obtain
\begin{equation*}
	\partial_t\dot{\Delta}_jS+\alpha_1\mathbf{u}\cdot\nabla\dot{\Delta}_jS+[\dot{\Delta}_j, \mathbf{u}\cdot\nabla]S=\dot{\Delta}_jf_3,
\end{equation*}
and followed the standard energy argument to get
\begin{equation}\label{2.54}
	\begin{split}
		||S^h||_{\widetilde{L}^\infty_t(\dot{B}^{\frac{n}{2}+1}_{2,1})}\lesssim&||S_0^h||_{\dot{B}^{\frac{n}{2}+1}_{2,1}}
+\int_{0}^{t}||(f_3)^h||_{\dot{B}^{\frac{n}{2}+1}_{2,1}}d\tau+\int_{0}^{t}||{\rm div}u||_{L^\infty}||S^h||_{\dot{B}^{\frac{n}{2}+1}_{2,1}}d\tau\\
		&+\int_{0}^{t}\sum\limits_{j\ge j_0}2^{(\frac{n}{2}+1)j}||[\dot{\Delta}_j, \mathbf{u}\cdot\nabla]S||_{L^2}d\tau.
	\end{split}
\end{equation}
By the expression of $f_3$, we have
\begin{equation}\label{2.55}
	||(f_3)^h||_{\dot{B}^{\frac{n}{2}+1}_{2,1}}\lesssim||(\Gamma(\alpha_1\mathbf{u}))^h||_{\dot{B}^{\frac{n}{2}+1}_{2,1}}
+||(I(P)\Gamma(\alpha_1\mathbf{u}))^h||_{\dot{B}^{\frac{n}{2}+1}_{2,1}}.
\end{equation}
On the basis of Bony decomposition and Proposition \ref{A.11}, we can bound $||(\Gamma(\alpha_1\mathbf{u}))^h||_{\dot{B}^{\frac{n}{2}+1}_{2,1}}$ as
\begin{equation}\label{2.56}
	\begin{split}
		||(\Gamma(\alpha_1\mathbf{u}))^h||_{\dot{B}^{\frac{n}{2}+1}_{2,1}}\lesssim&||(\nabla \mathbf{u})^2||_{\dot{B}^{\frac{n}{2}+1}_{2,1}}
		\lesssim||T_{\nabla \mathbf{u}}\nabla \mathbf{u}+R(\nabla \mathbf{u},\nabla \mathbf{u})||_{\dot{B}^{\frac{n}{2}+1}_{2,1}}\\
		\lesssim&||\nabla \mathbf{u}||_{\dot{B}^{\frac{n}{p}-1}_{p,1}}||\nabla \mathbf{u}||_{\dot{B}^{\frac{n}{p}+2}_{p,1}}
		\lesssim||\mathbf{u}||_{\dot{B}^{\frac{n}{p}}_{p,1}}||\mathbf{u}||_{\dot{B}^{\frac{n}{p}+3}_{p,1}}\\
		\lesssim&(||\mathbf{u}^l||_{\dot{B}^{\frac{n}{2}-1}_{2,1}}
+||\mathbf{u}^h||_{\dot{B}^{\frac{n}{p}+1}_{p,1}})(||\mathbf{u}^l||_{\dot{B}^{\frac{n}{2}+1}_{2,1}}+||\mathbf{u}^h||_{\dot{B}^{\frac{n}{p}+3}_{p,1}})\\
		\lesssim&\mathscr{E}_\infty(t)\mathscr{E}_1(t).
	\end{split}
\end{equation}
For $||(I(P)\Gamma(\alpha_1\mathbf{u}))^h||_{\dot{B}^{\frac{n}{2}+1}_{2,1}}$, Bony decomposition gives
\begin{equation*}
	I(P)\Gamma(\alpha_1\mathbf{u})=T_{I(P)}\Gamma(\alpha_1\mathbf{u})+R(I(P),\Gamma(\alpha_1\mathbf{u}))+T_{\Gamma(\alpha_1\mathbf{u})}I(P).
\end{equation*}
Then, we can get from Proposition \ref{A.11} and Proposition \ref{A.13} that
\begin{equation}\label{2.57}
	\begin{split}
		&||T_{I(P)}\Gamma(\alpha_1\mathbf{u})+R(I(P),\Gamma(\alpha_1\mathbf{u}))||_{\dot{B}^{\frac{n}{2}+1}_{2,1}}\\
		\lesssim&||I(P)||_{\dot{B}^{\frac{n}{p}-1}_{p,1}}||\Gamma(\alpha_1\mathbf{u})||_{\dot{B}^{\frac{n}{p}+2}_{p,1}}
		\lesssim||P||_{\dot{B}^{\frac{n}{p}-1}_{p,1}}\mathscr{E}_\infty(t)\mathscr{E}_1(t)\\
		\lesssim&(||P||^l_{\dot{B}^{\frac{n}{2}-1}_{2,1}}+||P||^h_{\dot{B}^{\frac{n}{p}+2}_{p,1}})\mathscr{E}_\infty(t)\mathscr{E}_1(t)\\
		\lesssim&(\mathscr{E}_\infty(t))^2\mathscr{E}_1(t),
	\end{split}
\end{equation}
where the term $||\Gamma(\alpha_1\mathbf{u})||_{\dot{B}^{\frac{n}{p}+2}_{p,1}}$ can be bounded as
\begin{equation}\label{2.58}
	\begin{split}
		||\Gamma(\alpha_1\mathbf{u})||_{\dot{B}^{\frac{n}{p}+2}_{p,1}}\lesssim&||(\nabla \mathbf{u})^2||_{\dot{B}^{\frac{n}{p}+2}_{p,1}}
		\lesssim||\nabla \mathbf{u}||_{L^\infty}||\nabla \mathbf{u}||_{\dot{B}^{\frac{n}{p}+2}_{p,1}}\\
		\lesssim&||\nabla \mathbf{u}||_{\dot{B}^{\frac{n}{p}}_{p,1}}||\nabla \mathbf{u}||_{\dot{B}^{\frac{n}{p}+2}_{p,1}}
		\lesssim||\mathbf{u}||_{\dot{B}^{\frac{n}{p}+1}_{p,1}}||\mathbf{u}||_{\dot{B}^{\frac{n}{p}+3}_{p,1}}\\
		\lesssim&(||\mathbf{u}^l||_{\dot{B}^{\frac{n}{2}-1}_{2,1}}+||\mathbf{u}^h||_{\dot{B}^{\frac{n}{p}+1}_{p,1}})(||\mathbf{u}^l||_{\dot{B}^{\frac{n}{2}+1}_{2,1}}||\mathbf{u}^h||_{\dot{B}^{\frac{n}{p}+3}_{p,1}})\\
		\lesssim&\mathscr{E}_\infty(t)\mathscr{E}_1(t).
	\end{split}
\end{equation}
And we also have
\begin{equation}\label{2.59}
	\begin{split}
		&||T_{\Gamma(\alpha_1\mathbf{u})}I(P)||_{\dot{B}^{\frac{n}{2}+1}_{2,1}}\\
\lesssim&||\Gamma(\alpha_1\mathbf{u})||_{\dot{B}^{\frac{n}{p}-1}_{p,1}}||I(P)||_{\dot{B}^{\frac{n}{p}+2}_{p,1}}
		\lesssim||(\nabla \mathbf{u})^2||_{\dot{B}^{\frac{n}{p}-1}_{p,1}}||P||_{\dot{B}^{\frac{n}{p}+2}_{p,1}}\\
		\lesssim&||\nabla \mathbf{u}||_{\dot{B}^{\frac{n}{p}}_{p,1}}||\nabla \mathbf{u}||_{\dot{B}^{\frac{n}{p}-1}_{p,1}}||P||_{\dot{B}^{\frac{n}{p}+2}_{p,1}}
		\lesssim||\mathbf{u}||_{\dot{B}^{\frac{n}{p}+1}_{p,1}}||\mathbf{u}||_{\dot{B}^{\frac{n}{p}}_{p,1}}||P||_{\dot{B}^{\frac{n}{p}+2}_{p,1}}\\
		\lesssim&(||\mathbf{u}^l||_{\dot{B}^{\frac{n}{2}+1}_{2,1}}+||\mathbf{u}^h||_{\dot{B}^{\frac{n}{p}+3}_{p,1}})
(||\mathbf{u}^l||_{\dot{B}^{\frac{n}{2}-1}_{2,1}}||\mathbf{u}^h||_{\dot{B}^{\frac{n}{p}+1}_{p,1}})\\
		&\times(||P^l||_{\dot{B}^{\frac{n}{2}-1}_{2,1}}+||P^h||_{\dot{B}^{\frac{n}{p}+2}_{p,1}})\\
		\lesssim&(\mathscr{E}_\infty(t))^2\mathscr{E}_1(t).
	\end{split}
\end{equation}
Combining \eqref{2.57} with \eqref{2.59}, we get
\begin{equation}\label{2.60}
	||(I(P)\Gamma(\alpha_1\mathbf{u}))^h||_{\dot{B}^{\frac{n}{2}+1}_{2,1}}\lesssim(\mathscr{E}_\infty(t))^2\mathscr{E}_1(t).
\end{equation}
Plugging \eqref{2.56} and \eqref{2.60} into \eqref{2.55}, we can acquire
\begin{equation}\label{2.61}
	||(f_3)^h||_{\dot{B}^{\frac{n}{2}+1}_{2,1}}\lesssim(\mathscr{E}_\infty(t))^2\mathscr{E}_1(t).
\end{equation}
For the remaining two terms in \eqref{2.54}, according to \eqref{a.17} of Lemma \ref{A.12}, it holds that
\begin{equation}\label{2.62}
	\begin{split}
		&||{\rm div}\mathbf{u}||_{L^\infty}||S^h||_{\dot{B}^{\frac{n}{2}+1}_{2,1}}
		+\sum\limits_{j\ge j_0}2^{(\frac{n}{2}+1)j}||[\dot{\Delta}_j, \mathbf{u}\cdot\nabla]S||_{L^2}\\
		\lesssim&||{\rm div}\mathbf{u}||_{\dot{B}^{\frac{n}{p}}_{p,1}}||S^h||_{\dot{B}^{\frac{n}{2}+1}_{2,1}}+||\nabla \mathbf{u}||_{\dot{B}^{\frac{n}{p}}_{p,1}}||S||_{\dot{B}^{\frac{n}{2}+1}_{2,1}}\\
		\lesssim&||\mathbf{u}||_{\dot{B}^{\frac{n}{p}+1}_{p,1}}||S^h||_{\dot{B}^{\frac{n}{2}+1}_{2,1}}
+||\mathbf{u}||_{\dot{B}^{\frac{n}{p}+1}_{p,1}}||S||_{\dot{B}^{\frac{n}{2}+1}_{2,1}}\\
		\lesssim&(||\mathbf{u}^l||_{\dot{B}^{\frac{n}{2}+1}_{2,1}}+||\mathbf{u}^h||_{\dot{B}^{\frac{n}{p}+3}_{p,1}})
(||S^l||_{\dot{B}^{\frac{n}{2}-1}_{2,1}}+||S^h||_{\dot{B}^{\frac{n}{2}+1}_{2,1}})\\
		\lesssim&\mathscr{E}_\infty(t)\mathscr{E}_1(t).
	\end{split}
\end{equation}
Finally, combining \eqref{2.61} with \eqref{2.62}, we can get from \eqref{2.54} that
\begin{equation}\label{2.63}
	||S^h||_{\widetilde{L}^\infty_t(\dot{B}^{\frac{n}{2}+1}_{p,1})}\lesssim||S_0^h||_{\dot{B}^{\frac{n}{2}+1}_{2,1}}+\int_{0}^{t}(1+\mathscr{E}_\infty(\tau)+(\mathscr{E}_\infty(\tau))^2)\mathscr{E}_\infty(\tau)\mathscr{E}_1(\tau)d\tau.
\end{equation}

\subsection{The estimate of $(\Lambda P,\mathbb{Q}\mathbf{u})$ in the high frequency.}
\quad\quad In this subsection, we are going to estimate the term $(\Lambda P,\mathbb{Q}\mathbf{u})$ in the high frequency. On the one hand, in order to get the results we want, we use a quantity from \cite{chen4,haspot} that
\begin{equation}\label{2.64}
	G:=\mathbb{Q}\mathbf{u}-\frac{\alpha_2}{\alpha_3+\alpha_4}\Delta^{-1}\nabla P,
\end{equation}
and then use the second equation of \eqref{2.1} to get
\begin{equation}\label{2.65}
	\partial_tG-(\alpha_3+\alpha_4)\Delta G=\mathbb{Q}f_2-\frac{\alpha_2}{\alpha_3+\alpha_4}\Delta^{-1}\nabla f_1+\frac{\alpha_2^2}{\alpha_3+\alpha_4}G+\frac{\alpha_2^3}{(\alpha_3+\alpha_4)^2}\Delta^{-1}\nabla P.
\end{equation}
Thus, we can gain by a standard energy argument that
\begin{equation}\label{2.66}
	\begin{split}
		&||G^h||_{\widetilde{L}^\infty_{t}(\dot{B}^{\frac{n}{p}+1}_{p,1})}+||G^h||_{L^1_t(\dot{B}^{\frac{n}{p}+3}_{p,1})}\\
		\lesssim&||G^h_0||_{\dot{B}^{\frac{n}{p}+1}_{p,1}}+\int_{0}^{t}||(\mathbb{Q}f_2)^h||_{\dot{B}^{\frac{n}{p}+1}_{p,1}}d\tau
+\int_{0}^{t}||(\Delta^{-1}\nabla f_1)^h||_{\dot{B}^{\frac{n}{p}+1}_{p,1}}d\tau\\
		&+\int_{0}^{t}||G^h||_{\dot{B}^{\frac{n}{p}+1}_{p,1}}d\tau+\int_{0}^{t}||(\Delta^{-1}\nabla P)^h||_{\dot{B}^{\frac{n}{p}+1}_{p,1}}d\tau\\
		\lesssim&||G^h_0||_{\dot{B}^{\frac{n}{p}+1}_{p,1}}+\int_{0}^{t}||(f_2)^h||_{\dot{B}^{\frac{n}{p}+1}_{p,1}}d\tau+\int_{0}^{t}||(f_1)^h||_{\dot{B}^{\frac{n}{p}}_{p,1}}d\tau\\
		&+\int_{0}^{t}||G^h||_{\dot{B}^{\frac{n}{p}+1}_{p,1}}d\tau+\int_{0}^{t}||P^h||_{\dot{B}^{\frac{n}{p}}_{p,1}}d\tau.\\
	\end{split}
\end{equation}
On the other hand, applying the operator ${\rm div}$ to \eqref{2.64} and sum up with the first equation of \eqref{2.1}, we have
\begin{equation*}
	\partial_tP+\frac{\alpha_2^2}{\alpha_3+\alpha_4}P+\alpha_1\mathbf{u}\cdot\nabla P=\frac{C_v+R}{C_v^2}\Gamma(\alpha_1\mathbf{u})-\frac{C_v+R}{C_v}\alpha_1P{\rm div}\mathbf{u}-\alpha_2{\rm div}G.
\end{equation*}
Taking $\dot{\Delta}_j\Lambda$ to the above equation and using the commutator's argument to get
\begin{equation}\label{2.67}
	\begin{split}
		&\partial_t\dot{\Delta}_j\Lambda P+\frac{\alpha_2^2}{\alpha_3+\alpha_4}\dot{\Delta}_j\Lambda P+\alpha_1\mathbf{u}\cdot\nabla \dot{\Delta}_j\Lambda P+\alpha_1[\dot{\Delta}_j, \mathbf{u}\cdot\nabla]\Lambda P\\
		=&-\alpha_1\dot{\Delta}_j(\Lambda \mathbf{u})\cdot\nabla P+\frac{C_v+R}{C_v^2}\dot{\Delta}_j\Lambda\Gamma(\alpha_1\mathbf{u})-\frac{C_v+R}{C_v}\alpha_1\dot{\Delta}_j\Lambda P{\rm div}\mathbf{u}-\alpha_2\dot{\Delta}_j\Lambda{\rm div}G.
	\end{split}
\end{equation}
Similarly, following a standard energy argument, we get
\begin{equation}\label{2.68}
	\begin{split}
		&||\Lambda P^h||_{\widetilde{L}^\infty_{t}(\dot{B}^{\frac{n}{p}+1}_{p,1})}+\frac{\alpha_2^2}{\alpha_3+\alpha_4}||\Lambda P^h||_{L^1_t(\dot{B}^{\frac{n}{p}+1}_{p,1})}\\
		\lesssim&||\Lambda P_0^h||_{\dot{B}^{\frac{n}{p}+1}_{p,1}}+\int_{0}^{t}||G^h||_{\dot{B}^{\frac{n}{p}+3}_{p,1}}d\tau\\
		& +\int_{0}^{t}||{\rm div}\mathbf{u}||_{L^\infty}||\Lambda P^h||_{\dot{B}^{\frac{n}{p}+1}_{p,1}}d\tau+\int_{0}^{t}\sum\limits_{j\ge j_0}2^{(\frac{n}{p}+1)j}||[\dot{\Delta}_j, \mathbf{u}\cdot\nabla]\Lambda P||_{L^2}d\tau\\
		& +\int_{0}^{t}||(P{\rm div}\mathbf{u})^h||_{\dot{B}^{\frac{n}{p}+2}_{p,1}}d\tau+\int_{0}^{t}||(\Lambda \mathbf{u}\cdot\nabla P)^h||_{\dot{B}^{\frac{n}{p}+1}_{p,1}}d\tau+\int_{0}^{t}||(\Gamma(\alpha_1\mathbf{u}))^h||_{\dot{B}^{\frac{n}{p}+2}_{p,1}}d\tau.
	\end{split}
\end{equation}
Because of the high frequency cut-off, we have
\begin{equation}\label{2.69}
	||G||^h_{L^1_t(\dot{B}^{\frac{n}{p}+1}_{p,1})}\lesssim2^{-2j_0}||G||^h_{L^1_t(\dot{B}^{\frac{n}{p}+3}_{p,1})}\ \ {\rm and}\ \ ||P||^h_{L^1_t(\dot{B}^{\frac{n}{p}}_{p,1})}\lesssim2^{-2j_0}||P||^h_{L^1_t(\dot{B}^{\frac{n}{p}+2}_{p,1})}.
\end{equation}
Then we can choose $j_0$ large enough and let $||G||^h_{L^1_t(\dot{B}^{\frac{n}{p}+1}_{p,1})}$, $||P||^h_{L^1_t(\dot{B}^{\frac{n}{p}}_{p,1})}$ are absorbed by the left hand side in \eqref{2.66} and \eqref{2.68}. Therefore, summing up \eqref{2.66}, \eqref{2.68}, we finally get
\begin{equation}\label{2.70}
	\begin{split}
		&||G^h||_{\widetilde{L}^\infty_{t}(\dot{B}^{\frac{n}{p}+1}_{p,1})}+||G^h||_{L^1_t(\dot{B}^{\frac{n}{p}+3}_{p,1})}+||\Lambda P^h||_{\widetilde{L}^\infty_{t}(\dot{B}^{\frac{n}{p}+1}_{p,1})}+\frac{\alpha_2^2}{\alpha_3+\alpha_4}||\Lambda P^h||_{L^1_t(\dot{B}^{\frac{n}{p}+1}_{p,1})}\\
		\lesssim&||G^h_0||_{\dot{B}^{\frac{n}{p}+1}_{p,1}}+||\Lambda P_0^h||_{\dot{B}^{\frac{n}{p}+1}_{p,1}}+\int_{0}^{t}||(f_1)^h||_{\dot{B}^{\frac{n}{p}}_{p,1}}d\tau+\int_{0}^{t}||(f_2)^h||_{\dot{B}^{\frac{n}{p}+1}_{p,1}}d\tau\\
		&+\int_{0}^{t}||{\rm div}\mathbf{u}||_{L^\infty}||\Lambda P^h||_{\dot{B}^{\frac{n}{p}+1}_{p,1}}d\tau+\int_{0}^{t}\sum\limits_{j\ge j_0}2^{(\frac{n}{p}+1)j}||[\dot{\Delta}_j, \mathbf{u}\cdot\nabla]\Lambda P||_{L^2}d\tau\\
		& +\int_{0}^{t}||(P{\rm div}\mathbf{u})^h||_{\dot{B}^{\frac{n}{p}+2}_{p,1}}d\tau+\int_{0}^{t}||(\Lambda \mathbf{u}\cdot\nabla P)^h||_{\dot{B}^{\frac{n}{p}+1}_{p,1}}d\tau+\int_{0}^{t}||(\Gamma(\alpha_1\mathbf{u}))^h||_{\dot{B}^{\frac{n}{p}+2}_{p,1}}d\tau.
	\end{split}
\end{equation}
For the terms of $||(f_1)^h||_{\dot{B}^{\frac{n}{p}}_{p,1}}$, we can obtain
\begin{equation}\label{2.71}
	||(f_1)^h||_{\dot{B}^{\frac{n}{p}}_{p,1}}\lesssim||(\mathbf{u}\cdot\nabla P)^h||_{\dot{B}^{\frac{n}{p}}_{p,1}}+||(P{\rm div}\mathbf{u})^h||_{\dot{B}^{\frac{n}{p}}_{p,1}}+||(\Gamma(\alpha_1\mathbf{u}))^h||_{\dot{B}^{\frac{n}{p}}_{p,1}}.
\end{equation}
Following \eqref{a.8} of Proposition \ref{A.7}, it holds
\begin{equation}\label{2.72}
	\begin{split}
		&||(\mathbf{u}\cdot\nabla P)^h||_{\dot{B}^{\frac{n}{p}}_{p,1}}+||(P{\rm div}\mathbf{u})^h||_{\dot{B}^{\frac{n}{p}}_{p,1}}\\
		\lesssim&||\mathbf{u}||_{\dot{B}^{\frac{n}{p}}_{p,1}}||\nabla P||_{\dot{B}^{\frac{n}{p}}_{p,1}}+||P||_{\dot{B}^{\frac{n}{p}}_{p,1}}||{\rm div}\mathbf{u}||_{\dot{B}^{\frac{n}{p}}_{p,1}}\\
		\lesssim&||\mathbf{u}||_{\dot{B}^{\frac{n}{p}}_{p,1}}||P||_{\dot{B}^{\frac{n}{p}+1}_{p,1}}
+||P||_{\dot{B}^{\frac{n}{p}}_{p,1}}||\mathbf{u}||_{\dot{B}^{\frac{n}{p}+1}_{p,1}}\\
		\lesssim&(||\mathbf{u}^l||_{\dot{B}^{\frac{n}{2}-1}_{2,1}}+||\mathbf{u}^h||_{\dot{B}^{\frac{n}{p}+1}_{p,1}})(||P^l||_{\dot{B}^{\frac{n}{2}+1}_{2,1}}
+||P^h||_{\dot{B}^{\frac{n}{p}+2}_{p,1}})\\
		&+(||P^l||_{\dot{B}^{\frac{n}{2}-1}_{2,1}}+||P^h||_{\dot{B}^{\frac{n}{p}+2}_{p,1}})(||\mathbf{u}^l||_{\dot{B}^{\frac{n}{2}+1}_{2,1}}
+||\mathbf{u}^h||_{\dot{B}^{\frac{n}{p}+3}_{p,1}})\\
		\lesssim&\mathscr{E}_\infty(t)\mathscr{E}_1(t),
	\end{split}
\end{equation}
and
\begin{equation}\label{2.73}
	\begin{split}
		||(\Gamma(\alpha_1\mathbf{u}))^h||_{\dot{B}^{\frac{n}{p}}_{p,1}}\lesssim&||(\nabla \mathbf{u})^2||_{\dot{B}^{\frac{n}{p}}_{p,1}}\\
		\lesssim&||\nabla \mathbf{u}||_{\dot{B}^{\frac{n}{p}}_{p,1}}||\nabla \mathbf{u}||_{\dot{B}^{\frac{n}{p}}_{p,1}}
		\lesssim||\mathbf{u}||_{\dot{B}^{\frac{n}{p}+1}_{p,1}}||\mathbf{u}||_{\dot{B}^{\frac{n}{p}+1}_{p,1}}\\
		\lesssim&(||\mathbf{u}^l||_{\dot{B}^{\frac{n}{2}-1}_{2,1}}+||\mathbf{u}^h||_{\dot{B}^{\frac{n}{p}+1}_{p,1}})
(||\mathbf{u}^l||_{\dot{B}^{\frac{n}{2}+1}_{2,1}}+||\mathbf{u}^h||_{\dot{B}^{\frac{n}{p}+3}_{p,1}})\\
		\lesssim&\mathscr{E}_\infty(t)\mathscr{E}_1(t).
	\end{split}
\end{equation}
So, plugging \eqref{2.72}, \eqref{2.73} into \eqref{2.71}, then $||(f_1)^h||_{\dot{B}^{\frac{n}{p}}_{p,1}}$ can be bounded as
\begin{equation}\label{2.74}
	||(f_1)^h||_{\dot{B}^{\frac{n}{p}}_{p,1}}\lesssim\mathscr{E}_\infty(t)\mathscr{E}_1(t).
\end{equation}
We have already estimated the term $||(f_2)^h||_{\dot{B}^{\frac{n}{p}+1}_{p,1}}$ completely in \eqref{2.52}. Now, by virtue of \eqref{a.17} of Lemma \ref{A.12}, we can get
\begin{equation}\label{2.75}
	\begin{split}
		&||{\rm div}\mathbf{u}||_{L^\infty}||\Lambda P^h||_{\dot{B}^{\frac{n}{p}+1}_{p,1}}+\sum\limits_{j\ge j_0}2^{(\frac{n}{p}+1)j}||[\dot{\Delta}_j, \mathbf{u}\cdot\nabla]\Lambda P||_{L^2}\\
		\lesssim&||{\rm div}\mathbf{u}||_{\dot{B}^{\frac{n}{p}}_{p,1}}||P^h||_{\dot{B}^{\frac{n}{p}+2}_{p,1}}+||\nabla \mathbf{u}||_{\dot{B}^{\frac{n}{p}}_{p,1}}||\Lambda P||_{\dot{B}^{\frac{n}{p}+1}_{p,1}}\\
		\lesssim&||\mathbf{u}||_{\dot{B}^{\frac{n}{p}+1}_{p,1}}||P^h||_{\dot{B}^{\frac{n}{p}+2}_{p,1}}
+||\mathbf{u}||_{\dot{B}^{\frac{n}{p}+1}_{p,1}}||P||_{\dot{B}^{\frac{n}{p}+2}_{p,1}}\\
		\lesssim&(||\mathbf{u}^l||_{\dot{B}^{\frac{n}{2}-1}_{2,1}}+||\mathbf{u}^h||_{\dot{B}^{\frac{n}{p}+1}_{p,1}})
(||P^l||_{\dot{B}^{\frac{n}{2}+1}_{2,1}}+||P^h||_{\dot{B}^{\frac{n}{p}+2}_{p,1}})\\
		\lesssim&\mathscr{E}_\infty(t)\mathscr{E}_1(t).
	\end{split}
\end{equation}
From \eqref{a.7} of Proposition \ref{A.7}, $||(P{\rm div}\mathbf{u})^h||_{\dot{B}^{\frac{n}{p}+2}_{p,1}}$ and $||(\Lambda \mathbf{u}\cdot\nabla P)^h||_{\dot{B}^{\frac{n}{p}+1}_{p,1}}$ can respectively be bounded as
\begin{equation}\label{2.76}
	\begin{split}
		||(P{\rm div}\mathbf{u})^h||_{\dot{B}^{\frac{n}{p}+2}_{p,1}}\lesssim&||P||_{L^\infty}||{\rm div}\mathbf{u}||_{\dot{B}^{\frac{n}{p}+2}_{p,1}}+||{\rm div}\mathbf{u}||_{L^\infty}||P||_{\dot{B}^{\frac{n}{p}+2}_{p,1}}\\
		\lesssim&||P||_{\dot{B}^{\frac{n}{p}}_{p,1}}||{\rm div}\mathbf{u}||_{\dot{B}^{\frac{n}{p}+2}_{p,1}}+||{\rm div}\mathbf{u}||_{\dot{B}^{\frac{n}{p}}_{p,1}}||P||_{\dot{B}^{\frac{n}{p}+2}_{p,1}}\\
		\lesssim&||P||_{\dot{B}^{\frac{n}{p}}_{p,1}}||\mathbf{u}||_{\dot{B}^{\frac{n}{p}+3}_{p,1}}
+||\mathbf{u}||_{\dot{B}^{\frac{n}{p}+1}_{p,1}}||P||_{\dot{B}^{\frac{n}{p}+2}_{p,1}}\\
		\lesssim&(||P^l||_{\dot{B}^{\frac{n}{2}-1}_{2,1}}
+||P^h||_{\dot{B}^{\frac{n}{p}+2}_{p,1}})(||\mathbf{u}^l||_{\dot{B}^{\frac{n}{2}+1}_{2,1}}+||\mathbf{u}^h||_{\dot{B}^{\frac{n}{p}+3}_{p,1}})\\
		\lesssim&\mathscr{E}_\infty(t)\mathscr{E}_1(t),
	\end{split}
\end{equation}
and
\begin{equation}\label{2.77}
	\begin{split}
		||(\Lambda \mathbf{u}\cdot\nabla P)^h||_{\dot{B}^{\frac{n}{p}+1}_{p,1}}\lesssim&||\Lambda \mathbf{u}||_{L^\infty}||\nabla P||_{\dot{B}^{\frac{n}{p}+1}_{p,1}}+||\nabla P||_{L^\infty}||\Lambda \mathbf{u}||_{\dot{B}^{\frac{n}{p}+1}_{p,1}}\\
		\lesssim&||\Lambda \mathbf{u}||_{\dot{B}^{\frac{n}{p}}_{p,1}}||\nabla P||_{\dot{B}^{\frac{n}{p}+1}_{p,1}}+||\nabla P||_{\dot{B}^{\frac{n}{p}}_{p,1}}||\Lambda \mathbf{u}||_{\dot{B}^{\frac{n}{p}+1}_{p,1}}\\
		\lesssim&||\mathbf{u}||_{\dot{B}^{\frac{n}{p}+1}_{p,1}}||P||_{\dot{B}^{\frac{n}{p}+2}_{p,1}}
+||P||_{\dot{B}^{\frac{n}{p}+1}_{p,1}}||\mathbf{u}||_{\dot{B}^{\frac{n}{p}+2}_{p,1}}\\
		\lesssim&(||\mathbf{u}^l||_{\dot{B}^{\frac{n}{2}+1}_{2,1}}+||\mathbf{u}^h||_{\dot{B}^{\frac{n}{p}+3}_{p,1}})
(||P^l||_{\dot{B}^{\frac{n}{2}-1}_{2,1}}+||P^h||_{\dot{B}^{\frac{n}{p}+2}_{p,1}})\\
		\lesssim&\mathscr{E}_\infty(t)\mathscr{E}_1(t).
	\end{split}
\end{equation}
For the term $||(\Gamma(\alpha_1\mathbf{u}))^h||_{\dot{B}^{\frac{n}{p}+2}_{p,1}}$, it has been estimated in \eqref{2.58}, so we have
\begin{equation}\label{2.78}
	||(\Gamma(\alpha_1\mathbf{u}))^h||_{\dot{B}^{\frac{n}{p}+2}_{p,1}}\lesssim\mathscr{E}_\infty(t)\mathscr{E}_1(t).
\end{equation}
Based on all the above inequalities, we have
\begin{equation}\label{2.79}
	\begin{split}
		&||G^h||_{\widetilde{L}^\infty_{t}(\dot{B}^{\frac{n}{p}+1}_{p,1})}+||G^h||_{L^1_t(\dot{B}^{\frac{n}{p}+3}_{p,1})}+||\Lambda P^h||_{\widetilde{L}^\infty_{t}(\dot{B}^{\frac{n}{p}+1}_{p,1})}+\frac{\alpha_2^2}{\alpha_3+\alpha_4}||\Lambda P^h||_{L^1_t(\dot{B}^{\frac{n}{p}+1}_{p,1})}\\
		\lesssim&||G^h_0||_{\dot{B}^{\frac{n}{p}+1}_{p,1}}+||\Lambda P_0^h||_{\dot{B}^{\frac{n}{p}+1}_{p,1}}
+\int_{0}^{t}(1+\mathscr{E}_\infty(\tau)+(\mathscr{E}_\infty(\tau))^2)\mathscr{E}_\infty(\tau)\mathscr{E}_1(\tau)d\tau.
	\end{split}
\end{equation}
Owing to $G:=\mathbb{Q}\mathbf{u}-\frac{\alpha_2}{\alpha_3+\alpha_4}\Delta^{-1}\nabla P$, we can simply get
\begin{equation}\label{2.80}
	\begin{split}
		||\mathbb{Q}\mathbf{u}^h||_{\widetilde{L}^\infty_{t}(\dot{B}^{\frac{n}{p}+1}_{p,1})}
\lesssim&||G^h||_{\widetilde{L}^\infty_{t}(\dot{B}^{\frac{n}{p}+1}_{p,1})}+||\Delta^{-1}\nabla P^h||_{\widetilde{L}^\infty_{t}(\dot{B}^{\frac{n}{p}+1}_{p,1})}\\
		\lesssim&||G^h||_{\widetilde{L}^\infty_{t}(\dot{B}^{\frac{n}{p}+1}_{p,1})}+||\Lambda P^h||_{\widetilde{L}^\infty_{t}(\dot{B}^{\frac{n}{p}+1}_{p,1})},
	\end{split}
\end{equation}
\begin{equation}\label{2.81}
	\begin{split}
		||\mathbb{Q}\mathbf{u}^h||_{L^1_t(\dot{B}^{\frac{n}{p}+3}_{p,1})}
\lesssim&||G^h||_{L^1_t(\dot{B}^{\frac{n}{p}+3}_{p,1})}+||\Delta^{-1}\nabla P^h||_{L^1_t(\dot{B}^{\frac{n}{p}+3}_{p,1})}\\
		\lesssim&||G^h||_{L^1_t(\dot{B}^{\frac{n}{p}+3}_{p,1})}+||\Lambda P^h||_{L^1_t(\dot{B}^{\frac{n}{p}+1}_{p,1})},
	\end{split}
\end{equation}
and
\begin{equation}\label{2.82}
	\begin{split}
		||G_0^h||_{\dot{B}^{\frac{n}{p}+1}_{p,1}}\lesssim&||\mathbb{Q}\mathbf{u}_0^h||_{\dot{B}^{\frac{n}{p}+1}_{p,1}}+||\Delta^{-1}\nabla P_0^h||_{\dot{B}^{\frac{n}{p}+1}_{p,1}}
		\lesssim||\mathbb{Q}\mathbf{u}_0^h||_{\dot{B}^{\frac{n}{p}+1}_{p,1}}+||\Lambda P_0^h||_{\dot{B}^{\frac{n}{p}+1}_{p,1}}.
	\end{split}
\end{equation}
Plugging \eqref{2.80}-\eqref{2.82} into \eqref{2.79}, we finally arrive at
\begin{equation}\label{2.83}
	\begin{split}	&||\mathbb{Q}\mathbf{u}^h||_{\widetilde{L}^\infty_{t}(\dot{B}^{\frac{n}{p}+1}_{p,1})}+||\mathbb{Q}\mathbf{u}^h||_{L^1_t(\dot{B}^{\frac{n}{p}+3}_{p,1})}+||\Lambda P^h||_{\widetilde{L}^\infty_{t}(\dot{B}^{\frac{n}{p}+1}_{p,1})}+\frac{\alpha_2^2}{\alpha_3+\alpha_4}||\Lambda P^h||_{L^1_t(\dot{B}^{\frac{n}{p}+1}_{p,1})}\\
		\lesssim&||\mathbb{Q}\mathbf{u}_0^h||_{\dot{B}^{\frac{n}{p}+1}_{p,1}}+||\Lambda P_0^h||_{\dot{B}^{\frac{n}{p}+1}_{p,1}}+\int_{0}^{t}(1+\mathscr{E}_\infty(\tau)+(\mathscr{E}_\infty(\tau))^2)\mathscr{E}_\infty(\tau)\mathscr{E}_1(\tau)d\tau.
	\end{split}
\end{equation}

\vspace{3mm}
			
\noindent{\textbf{The proof of Theorem \ref{Theorem 1.1}}}:
First, define
\begin{equation}
\begin{array}{rl}
&\mathcal{X}(t)\triangleq\sup\limits_{t>0}\mathcal{E}_\infty(t)+\int_0^t\mathcal{E}_1(\tau)d\tau,\\[2mm]
&\displaystyle\mathcal{X}(0)\triangleq\|(P_0^l,\mathbf{u}_0^l,S_0^l)\|_{\dot{B}_{2,1}^{\frac{n}{2}-1}}+\|S_0^h\|_{\dot{B}_{2,1}^{\frac{n}{2}+1}}+\|(\Lambda P_0^h,\mathbf{u}_0^h)\|_{\dot{B}_{p,1}^{\frac{n}{p}+1}}.
\end{array}
\end{equation}
Then, by virtue of the estimates in Sections 2.1-2.6, one knows that there exists a constant $C>0$ independent of $t$ such that
\begin{equation}\label{2.85}
\mathcal{X}(t)\leq C\mathcal{X}(0)+C(\mathcal{X}(t)^2+\mathcal{X}(t)^3).
\end{equation}

On the other hand, the local existence can be achieved via standard argument as in Haspot \cite{haspot}. Using the setting of initial data in Theorem \ref{Theorem 1.1} and the local existence, one can get for some $C>0$ that
\begin{equation}\label{2.86}
\mathcal{X}(t)\leq 2Cc_0,\ \ \forall t\in[0,T].
\end{equation}
Then, by taking a standard continuation argument, we can complete the proof of Theorem \ref{Theorem 1.1}.

\section{Decay rate}

\quad\quad In this section, we shall derive the decay rate of the solution constructed in Section 2.

From Section 2, we can get the following inequality:
\begin{equation}\label{4.28}
\begin{array}{ll}
&\frac{d}{dt}\Big(\|(P,\mathbf{u})^l\|_{\dot B_{2,1}^{\frac{n}{2}-1}}+\|(\Lambda P,\mathbf{u})\|_{\dot B_{p,1}^{\frac{n}{p}+1}}^h\Big)+\Big(\|(P,\mathbf{u})^l\|_{\dot B_{2,1}^{\frac{n}{2}+1}}+\|\Lambda P\|_{\dot B_{p,1}^{\frac{n}{p}+1}}^h+\|\mathbf{u}\|_{\dot B_{p,1}^{\frac{n}{p}+3}}^h\Big)\\[2mm]
&\quad\lesssim \Big(1+\|(P,\mathbf{u},S)\|_{\dot B_{2,1}^{\frac{n}{2}-1}}^l+\|\Lambda P\|_{\dot B_{p,1}^{\frac{n}{p}+1}}^h+\|S\|_{\dot B_{2,1}^{\frac{n}{2}+1}}^h+\|\mathbf{u}\|_{\dot B_{p,1}^{\frac{n}{p}+1}}^h\Big)\\[2mm]
&\quad\quad\times \Big(\|(P,\mathbf{u},S)\|_{\dot B_{2,1}^{\frac{n}{2}-1}}^l+\|\Lambda P\|_{\dot B_{p,1}^{\frac{n}{p}+1}}^h+\|S\|_{\dot B_{2,1}^{\frac{n}{2}+1}}^h+\|\mathbf{u}\|_{\dot B_{p,1}^{\frac{n}{p}+1}}^h\Big)\\[2mm]
&\quad\quad\times \Big(\|(P,\mathbf{u})\|_{\dot B_{2,1}^{\frac{n}{2}+1}}^l+\|P\|_{\dot B_{p,1}^{\frac{n}{p}+2}}^h+\|\mathbf{u}\|_{\dot B_{p,1}^{\frac{n}{p}+1}}^h\Big).
\end{array}
\end{equation}
By the proof of the global existence of Theorem \ref{Theorem 1.1}, the following estimate holds:
\begin{equation}\label{4.29}
\|(P,\mathbf{u},S)\|_{\dot B_{2,1}^{\frac{n}{2}-1}}^l+\|\Lambda P\|_{\dot B_{p,1}^{\frac{n}{p}+1}}^h+\|S\|_{\dot B_{2,1}^{\frac{n}{2}+1}}^h+\|\mathbf{u}\|_{\dot B_{p,1}^{\frac{n}{p}+1}}^h\leq C_0,
\end{equation}
from which we can infer from (\ref{4.28}) that there exists a constant $\bar c>0$ such that
\begin{equation}\label{4.30}
\frac{d}{dt}\Big(\|(P,\mathbf{u})^l\|_{\dot B_{2,1}^{\frac{n}{2}-1}}+\|(\Lambda P,\mathbf{u})\|_{\dot B_{p,1}^{\frac{n}{p}+1}}^h\Big)+\bar{c}\Big(\|(P,\mathbf{u})^l\|_{\dot B_{2,1}^{\frac{n}{2}+1}}+\|\Lambda P\|_{\dot B_{p,1}^{\frac{n}{p}+1}}^h+\|\mathbf{u}\|_{\dot B_{p,1}^{\frac{n}{p}+3}}^h\Big)\leq 0.
\end{equation}
With (\ref{4.30}) in hand, one can get a Lyapunov-type differential inequality, which relies heavily on an interpolation inequality. Before using the interpolation inequality, it requires the uniform bound as follows:
\begin{equation}\label{4.31}
\|(P,\mathbf{u})\|_{\dot B_{2,\infty}^\sigma}^l\leq C,\ for\ any\ \frac{n}{2}-\frac{2n}{p}\leq\sigma<\frac{n}{2}-1.
\end{equation}

\subsection{Propagation the regularity of the initial data with negative index}
\quad\quad In this subsection, we shall derive the following key proposition.
\begin{proposition}\label{proposition 4.1}
Let $(P,\mathbf{u},S)$ be the solutions constructed in Section $2$. For any $\frac{n}{2}-\frac{2n}{p}\leq\sigma<\frac{n}{2}-1$ and $(P_0^l,\mathbf{u}_0^l,S_0^l)\in\dot{B}_{2,\infty}^\sigma$, then there exists a constant $C_0>0$ depends on the norm of the initial data such that for all $t\geq 0$,
\begin{equation}\label{4.1}
\|(P,\mathbf{u},S)(t,\cdot)\|_{\dot{B}_{2,\infty}^\sigma}^l\leq C_0.
\end{equation}
\begin{proof}
For the first two equations of (\ref{2.1}), we get that
\begin{equation}\label{4.2}
\|(P,\mathbf{u})\|_{\dot{B}_{2,\infty}^\sigma}^l\lesssim\|(P_0,\mathbf{u}_0)\|_{\dot{B}_{2,\infty}^\sigma}^l+\int_0^t\|(f_1,f_2)\|_{\dot{B}_{2,\infty}^\sigma}^ld\tau.
\end{equation}
In what follows, we focus on the nonlinear norm $\|(f_1,f_2)\|_{\dot{B}_{2,\infty}^\sigma}^l$. By using Corollary \ref{A.9} and the decompositions $\mathbf{u}=\mathbf{u}^l+\mathbf{u}^h,\ P=P^l+P^h$, we find that
\begin{equation}\label{4.3}
\begin{array}{ll}
&\|\mathbf{u}\cdot\nabla P^l\|_{\dot{B}_{2,\infty}^\sigma}^l+\|P{\rm div}\mathbf{u}^l\|_{\dot{B}_{2,\infty}^\sigma}^l\\[3mm]
\lesssim&\|\mathbf{u}^l\cdot\nabla P^l\|_{\dot{B}_{2,\infty}^\sigma}+\|\mathbf{u}^h\cdot\nabla P^l\|_{\dot{B}_{2,\infty}^\sigma}+\|P^l{\rm div}\mathbf{u}^l\|_{\dot{B}_{2,\infty}^\sigma}+\|P^h{\rm div}\mathbf{u}^l\|_{\dot{B}_{2,\infty}^\sigma}\\[3mm]
\lesssim&\|\mathbf{u}^l\|_{\dot{B}_{2,\infty}^\sigma}\|\nabla P^l\|_{\dot{B}_{p,1}^{\frac{n}{p}}}+\|\mathbf{u}^h\|_{\dot{B}_{p,1}^{\frac{n}{p}}}\|\nabla P^l\|_{\dot{B}_{2,\infty}^\sigma}\\[3mm]
&+\|P^l\|_{\dot{B}_{2,\infty}^\sigma}\|{\rm div}\mathbf{u}^l\|_{\cdot{B}_{p,1}^{\frac{n}{p}}}+\|P^h\|_{\dot{B}_{p,1}^{\frac{n}{p}}}\|{\rm div}\mathbf{u}^l\|_{\dot{B}_{2,\infty}^\sigma}\\[3mm]
\lesssim&\|\mathbf{u}^l\|_{\dot{B}_{2,\infty}^\sigma}\|P^l\|_{\dot{B}_{p,1}^{\frac{n}{p}+1}}+\|\mathbf{u}^h\|_{\dot{B}_{p,1}^{\frac{n}{p}+3}}\|P^l\|_{\dot{B}_{2,\infty}^\sigma}\\[3mm]
&\quad\quad\quad+\|P^l\|_{\dot{B}_{2,\infty}^\sigma}\|\mathbf{u}^l\|_{\dot{B}_{p,1}^{\frac{n}{p}+1}}
+\|P^h\|_{\dot{B}_{p,1}^{\frac{n}{p}+2}}\|\mathbf{u}^l\|_{\dot{B}_{2,\infty}^\sigma}\\[3mm]
\lesssim&\mathcal{E}_1(t)\|(P,\mathbf{u})\|_{\dot{B}_{2,\infty}^\sigma}^l.
\end{array}
\end{equation}
According to Proposition \ref{A.10} with $s=\frac{n}{p}$, it holds that
\begin{equation}\label{4.0}
\|fg^h\|_{\dot B_{2,\infty}^{\frac{n}{2}-\frac{2n}{p}}}^l\lesssim(\|f\|_{\dot B_{p,1}^s}+\|f^l\|_{L^{p^\ast}})\|g^h\|_{\dot B_{p,1}^{-s}},
\end{equation}
where $\frac{1}{p}+\frac{1}{p^\ast}=\frac{1}{2}$.
By using (\ref{4.0}) and $\dot{B}_{2,1}^{\frac{n}{p}}\hookrightarrow L^{p^\ast}$, we have
\begin{equation}\label{4.4}
\begin{array}{ll}
&\|\mathbf{u}\cdot\nabla P^h\|_{\dot{B}_{2,\infty}^\sigma}^l+\|P{\rm div}\mathbf{u}^h\|_{\dot{B}_{2,\infty}^\sigma}^l\\[3mm]
&\quad\quad\lesssim\|\mathbf{u}\cdot\nabla P^h\|_{\dot{B}_{2,\infty}^{\frac{n}{2}-\frac{2n}{p}}}^l+\|P{\rm div}\mathbf{u}^h\|_{\dot{B}_{2,\infty}^{\frac{n}{2}-\frac{2n}{p}}}^l\\[3mm]
&\quad\quad\lesssim(\|\mathbf{u}\|_{\dot{B}_{p,1}^s}+\|\mathbf{u}^l\|_{L^{p^\ast}})\|\nabla P^h\|_{\dot{B}_{p,1}^{-s}}+(\|P\|_{\dot{B}_{p,1}^s}+\|P^l\|_{L^{p^\ast}})\|{\rm div} \mathbf{u}^h\|_{\dot{B}_{p,1}^{-s}}\\[3mm]
&\quad\quad\lesssim(\|\mathbf{u}\|_{\dot{B}_{p,1}^s}+\|\mathbf{u}^l\|_{\dot{B}_{2,1}^{\frac{n}{p}}})\|\nabla P^h\|_{\dot{B}_{p,1}^{-s}}+(\|P\|_{\dot{B}_{p,1}^s}+\|P^l\|_{\dot{B}_{2,1}^{\frac{n}{p}}})\|{\rm div} \mathbf{u}^h\|_{\dot{B}_{p,1}^{-s}}\\[3mm]
&\quad\quad\lesssim\mathcal{E}_\infty(t)\mathcal{E}_1(t).
\end{array}
\end{equation}
Similarly, we estimate the last term of $f_1$ as follows:
\begin{equation}\label{4.5}
\begin{array}{ll}
\|\frac{\Gamma(\alpha_1\mathbf{u})}{C_v}\|_{\dot{B}_{2,\infty}^\sigma}^l\!\!\!&\lesssim\|(\nabla \mathbf{u})^2\|_{\dot{B}_{2,\infty}^\sigma}^l\\[2mm]
\!\!\!&\lesssim\|\nabla \mathbf{u}\cdot\nabla \mathbf{u}^l\|_{\dot{B}_{2,\infty}^\sigma}^l+\|\nabla \mathbf{u}\cdot\nabla \mathbf{u}^h\|_{\dot{B}_{2,\infty}^\sigma}^l\\[2mm]
\!\!\!&\lesssim\|\nabla \mathbf{u}^l\|_{\dot{B}_{p,1}^{\frac{n}{p}}}\|\nabla \mathbf{u}^l\|_{\dot{B}_{2,\infty}^\sigma}+\|\nabla \mathbf{u}^h\|_{\dot{B}_{p,1}^{\frac{n}{p}}}\|\nabla \mathbf{u}^l\|_{\dot{B}_{2,\infty}^\sigma}+\|\nabla \mathbf{u}\cdot\nabla \mathbf{u}^h\|_{\dot{B}_{2,\infty}^{\frac{n}{2}-\frac{2n}{p}}}^l\\[2mm]
\!\!\!&\lesssim\mathcal{E}_1(t)\|\mathbf{u}^l\|_{\dot{B}_{2,\infty}^\sigma}+(\|\nabla \mathbf{u}\|_{\dot{B}_{p,1}^s}+\|\nabla \mathbf{u}^l\|_{L^{p^\ast}})\|\nabla \mathbf{u}^h\|_{\dot{B}_{p,1}^{-s}}\\[2mm]
&\lesssim\mathcal{E}_1(t)\|\mathbf{u}^l\|_{\dot{B}_{2,\infty}^\sigma}+(\|\nabla \mathbf{u}\|_{\dot{B}_{p,1}^s}+\|\nabla \mathbf{u}^l\|_{\dot{B}_{2,1}^{\frac{n}{p}}})\|\nabla \mathbf{u}^h\|_{\dot{B}_{p,1}^{-s}}\\[2mm]
&\lesssim\mathcal{E}_1(t)\|\mathbf{u}\|_{\dot{B}_{2,\infty}^\sigma}^l+\mathcal{E}_\infty(t)\mathcal{E}_1(t).
\end{array}
\end{equation}
We now turn to bound the terms in $f_2$. At first, we can use the same way as $\mathbf{u}\cdot\nabla P$ and $P{\rm div}\mathbf{u}$ to estimate $\mathbf{u}\cdot\nabla \mathbf{u}$ so that
\begin{equation}\label{4.6}
\begin{array}{ll}
\|\mathbf{u}\cdot\nabla \mathbf{u}\|_{\dot{B}_{2,\infty}^\sigma}^l&\lesssim\|\nabla \mathbf{u}^l\|_{\dot{B}_{p,1}^{\frac{n}{p}}}\|\mathbf{u}^l\|_{\dot{B}_{2,\infty}^\sigma}+\|\mathbf{u}^h\|_{\dot{B}_{p,1}^{\frac{n}{p}}}\|\nabla \mathbf{u}^l\|_{\dot{B}_{2,\infty}^\sigma}\\[2mm]
&\quad+(\|\mathbf{u}\|_{\dot{B}_{p,1}^s}+\|\mathbf{u}^l\|_{\dot{B}_{2,1}^{\frac{n}{p}}})\|\nabla \mathbf{u}^h\|_{\dot{B}_{p,1}^{-s}}\\[2mm]
&\lesssim\mathcal{E}_1(t)\|\mathbf{u}\|_{\dot{B}_{2,\infty}^\sigma}^l+\mathcal{E}_\infty(t)\mathcal{E}_1(t).
\end{array}
\end{equation}
Noting that
\begin{equation}\label{4.7}
\begin{array}{cc}
k(\rho)=k'(0)\rho+\rho\tilde{k}(\rho),\ \ \tilde{k}(0)=0,\\[2mm]
\frac{1}{\rho}-\frac{1}{\bar\rho}\sim o(1)(P+S),
\end{array}
\end{equation}
we can bound the term $\|k(\rho)\nabla\rho\|_{\dot{B}_{2,\infty}^\sigma}^l$ as
\begin{equation}\label{4.8}
\begin{array}{ll}
\|k(\rho)\nabla P\|_{\dot{B}_{2,\infty}^\sigma}^l&\lesssim\|k'(0)\rho\nabla P\|_{\dot{B}_{2,\infty}^\sigma}^l+\|\rho\tilde{k}(\rho)\nabla P\|_{\dot{B}_{2,\infty}^\sigma}^l\\[2mm]
&\lesssim\|k'(0)P\nabla P\|_{\dot{B}_{2,\infty}^\sigma}^l+\|k'(0)S\nabla P\|_{\dot{B}_{2,\infty}^\sigma}^l\\[2mm]
&\quad+\|P\tilde{k}(\rho)\nabla P\|_{\dot{B}_{2,\infty}^\sigma}^l+\|S\tilde{k}(\rho)\nabla P\|_{\dot{B}_{2,\infty}^\sigma}^l.
\end{array}
\end{equation}
To deal with the first term on the right side of (\ref{4.8}), we use Corollary \ref{A.9} and (\ref{4.0}) to get
\begin{equation}\label{4.9}
\begin{array}{ll}
\|k'(0)P\nabla P\|_{\dot{B}_{2,\infty}^\sigma}^l&\lesssim\|P^l\nabla P\|_{\dot{B}_{2,\infty}^\sigma}^l+\|P^h\nabla P\|_{\dot{B}_{2,\infty}^\sigma}^l\\[2mm]
&\lesssim\|P^l\|_{\dot{B}_{2,\infty}^\sigma}\|\nabla P\|_{\dot{B}_{p,1}^{\frac{n}{p}}}+(\|\nabla P\|_{\dot{B}_{p,1}^s}+\|\nabla P^l\|_{L^{p^\ast}})\|P^h\|_{\dot{B}_{p,1}^{-s}}\\[2mm]
&\lesssim\mathcal{E}_1(t)\|P^l\|_{\dot{B}_{2,\infty}^\sigma}+(\|\nabla P\|_{\dot{B}_{p,1}^s}+\|\nabla P^l\|_{\dot{B}_{2,1}^{\frac{n}{p}}})\|P^h\|_{\dot{B}_{p,1}^{-s}}\\[2mm]
&\lesssim\mathcal{E}_1(t)\|P^l\|_{\dot{B}_{2,\infty}^\sigma}+\mathcal{E}_\infty(t)\mathcal{E}_1(t).
\end{array}
\end{equation}
According to the embedding relation $\dot B_{2,1}^\sigma\hookrightarrow\dot B_{2,\infty}^\sigma$ and Proposition \ref{A.8}, one has
\begin{equation}\label{4.10}
\begin{array}{ll}
\|k'(0)S\nabla P\|_{\dot{B}_{2,\infty}^\sigma}^l\!\!\!&\lesssim\|S\nabla P^l\|_{\dot{B}_{2,\infty}^\sigma}^l+\|S\nabla P^h\|_{\dot{B}_{2,\infty}^\sigma}^l\\[2mm]
\!\!\!&\lesssim\|S\|_{\dot{B}_{2,\infty}^\sigma}\|\nabla P^l\|_{\dot{B}_{p,1}^{\frac{n}{p}}}+\|S\nabla P^h\|_{\dot{B}_{2,\infty}^{\sigma+\frac{n}{p}-\frac{n}{2}}}^l\\[2mm]
\!\!\!&\lesssim(\|S^l\|_{\dot{B}_{2,\infty}^\sigma}\!+\!\|S^h\|_{\dot{B}_{2,1}^{\frac{n}{2}}})\|P^l\|_{\dot{B}_{p,1}^{\frac{n}{p}+1}}\!+\!\|S\|_{\dot{B}_{2,\infty}^{\frac{n}{p}-\frac{n}{2}+\sigma+1}}+\|\nabla P^h\|_{\dot{B}_{p,1}^{\frac{n}{p}-1}}\\[2mm]
\!\!\!&\lesssim(\|S^l\|_{\dot{B}_{2,\infty}^\sigma}+\mathcal{E}_\infty(t))\mathcal{E}_1(t)+(\|S^l\|_{\dot{B}_{2,\infty}^\sigma}+\|S^h\|_{\dot{B}_{2,1}^{\frac{n}{2}+1}})\|P^h\|_{\dot{B}_{p,1}^{\frac{n}{p}+2}}\\[2mm]
\!\!\!&\lesssim\mathcal{E}_1(t)\|S^l\|_{\dot{B}_{2,\infty}^\sigma}+\mathcal{E}_\infty(t)\mathcal{E}_1(t).
\end{array}
\end{equation}
Similarly, by using embedding relation, Proposition \ref{A.8} and (\ref{4.0}), then $\|P\tilde{k}(\rho)\nabla P\|_{\dot{B}_{2,\infty}^\sigma}^l$ can be bounded as
\begin{equation}\label{4.11}
\begin{array}{ll}
\|P\tilde{k}(\rho)\nabla P\|_{\dot{B}_{2,\infty}^\sigma}^l&\lesssim\|P^l\tilde{k}(\rho)\nabla P\|_{\dot{B}_{2,\infty}^\sigma}^l+\|P^h\tilde{k}(\rho)\nabla P\|_{\dot{B}_{2,\infty}^\sigma}^l\\[2mm]
&\lesssim\|P^l\|_{\dot{B}_{2,\infty}^\sigma}\|\tilde{k}(\rho)\nabla P\|_{\dot{B}_{p,1}^{\frac{n}{p}}}\\[2mm]
&+(\|\tilde{k}(\rho)\nabla P\|_{\dot{B}_{p,1}^s}+\|(\tilde{k}(\rho)\nabla P)^l\|_{L^{p^\ast}})\|P^h\|_{\dot{B}_{p,1}^{-s}}\\[2mm]
&\lesssim\|P^l\|_{\dot{B}_{2,\infty}^\sigma}\|\tilde{k}(\rho)\|_{\dot{B}_{p,1}^{\frac{n}{p}}}\|\nabla P\|_{\dot{B}_{p,1}^{\frac{n}{p}}}+(\|\tilde{k}(\rho)\|_{\dot{B}_{p,1}^{\frac{n}{p}}}\|\nabla P\|_{\dot{B}_{p,1}^s}\\[2mm]
&\quad+\|\tilde{k}(\rho)\nabla P^l\|_{\dot{B}_{2,1}^{\frac{n}{p}}}+\|\tilde{k}(\rho)\nabla P^h\|_{\dot{B}_{p,1}^{\frac{2n}{p}-\frac{n}{2}}})\|P^h\|_{\dot{B}_{p,1}^{-s}}\\[2mm]
&\lesssim\|P^l\|_{\dot{B}_{2,\infty}^\sigma}\|\tilde{k}(\rho)\|_{\dot{B}_{p,1}^{\frac{n}{p}}}\|\nabla P\|_{\dot{B}_{p,1}^{\frac{n}{p}}}+[\|\tilde{k}(\rho)\|_{\dot{B}_{p,1}^{\frac{n}{p}}}\|\nabla P\|_{\dot{B}_{p,1}^s}\\[2mm]
&\quad+\|\tilde{k}(\rho)\|_{\dot{B}_{p,1}^{\frac{n}{p}}}(\|\nabla P^l\|_{\dot{B}_{2,1}^{\frac{n}{p}}}+\|\nabla P^h\|_{\dot{B}_{p,1}^{\frac{2n}{p}-\frac{n}{2}}})]\|P^h\|_{\dot{B}_{p,1}^{-s}}\\[2mm]
&\lesssim\mathcal{E}_\infty(t)\mathcal{E}_1(t)\|P^l\|_{\dot{B}_{2,\infty}^\sigma}+(\mathcal{E}_{\infty}(t))^2\mathcal{E}_1(t),
\end{array}
\end{equation}
and $\|S\tilde{k}(\rho)\nabla P\|_{\dot{B}_{2,\infty}^\sigma}^l$ can be bounded as
\begin{equation}\label{4.12}
\begin{array}{ll}
\!\!\!\!\!\!\!\|S\tilde{k}(\rho)\nabla P\|_{\dot{B}_{2,\infty}^\sigma}^l\!\!\!&\lesssim\|S\tilde{k}(\rho)\nabla P^l\|_{\dot{B}_{2,\infty}^\sigma}^l+\|S\tilde{k}(\rho)\nabla P^h\|_{\dot{B}_{2,\infty}^\sigma}^l\\[2mm]
\!\!\!&\lesssim\|S\tilde{k}(\rho)\|_{\dot{B}_{2,\infty}^\sigma}\|\nabla P^l\|_{\dot{B}_{p,1}^{\frac{n}{p}}}\!+\!(\|S\tilde{k}(\rho)\|_{\dot{B}_{p,1}^s}\!+\!\|(S\tilde{k}(\rho))^l\|_{L^{p^\ast}})\|\nabla P^h\|_{\dot{B}_{p,1}^{-s}}\\[2mm]
\!\!\!&\lesssim\|S\|_{\dot{B}_{2,\infty}^\sigma}\|\tilde{k}(\rho)\|_{\dot{B}_{p,1}^{\frac{n}{p}}}\|\nabla P^l\|_{\dot{B}_{p,1}^{\frac{n}{p}}}\\[2mm]
\!\!\!&\quad+(\|S\|_{\dot{B}_{p,1}^{\frac{n}{p}}}\|\tilde{k}(\rho)\|_{\dot{B}_{p,1}^s}+\|(S\tilde{k}(\rho))^l\|_{\dot{B}_{2,1}^{\frac{n}{p}}})\|\nabla P^h\|_{\dot{B}_{p,1}^{-s}}\\[2mm]
\!\!\!&\lesssim\|S\|_{\dot{B}_{2,\infty}^\sigma}\|\tilde{k}(\rho)\|_{\dot{B}_{p,1}^{\frac{n}{p}}}\|P^l\|_{\dot{B}_{p,1}^{\frac{n}{p}+1}}\\[2mm]
\!\!\!&\quad+(\|S\|_{\dot{B}_{p,1}^{\frac{n}{p}}}\|\tilde{k}(\rho)\|_{\dot{B}_{p,1}^s}+\|\tilde{k}(\rho)\|_{\dot{B}_{p,1}^{\frac{n}{p}}}\|S\|_{\dot{B}_{2,1}^{\frac{n}{p}}})\|\nabla P^h\|_{\dot{B}_{p,1}^{-s}}\\[2mm]
\!\!\!&\lesssim(\|S^l\|_{\dot{B}_{2,\infty}^\sigma}+\|S^h\|_{\dot{B}_{2,1}^{\frac{n}{2}+1}})\mathcal{E}_\infty(t)\mathcal{E}_1(t)+(\mathcal{E}_\infty(t))^2\mathcal{E}_1(t)\\[2mm]
\!\!\!&\lesssim\mathcal{E}_\infty(t)\mathcal{E}_1(t)\|S^l\|_{\dot{B}_{2,\infty}^\sigma}+(\mathcal{E}_{\infty}(t))^2\mathcal{E}_1(t).
\end{array}
\end{equation}
As a result, we get
\begin{equation}\label{4.13}
\begin{array}{ll}
\|k(\rho)\nabla P\|_{\dot{B}_{2,\infty}^\sigma}^l\lesssim(1+\mathcal{E}_\infty(t)\mathcal{E}_1(t))\|(P,S)\|_{\dot{B}_{2,\infty}^\sigma}^l+(1+\mathcal{E}_{\infty}(t))\mathcal{E}_{\infty}(t)\mathcal{E}_1(t).
\end{array}
\end{equation}
Then,  consider $\|k(\rho)(\mu\Delta \mathbf{u}+(\lambda+\mu)\nabla{\rm div}\mathbf{u}\|_{\dot{B}_{2,\infty}^\sigma}^l$. Similar to (\ref{4.8}), we obtain
\begin{equation}\label{4.14}
\begin{array}{ll}
\|k(\rho)(\mu\Delta \mathbf{u}+(\lambda+\mu)\nabla{\rm div}\mathbf{u}\|_{\dot{B}_{2,\infty}^\sigma}^l&\lesssim\|k(\rho)\nabla^2\mathbf{u}\|_{\dot{B}_{2,\infty}^\sigma}^l\\[2mm]
&\lesssim\|k'(0)P\nabla^2\mathbf{u}\|_{\dot{B}_{2,\infty}^\sigma}^l+\|k'(0)S\nabla^2\mathbf{u}\|_{\dot{B}_{2,\infty}^\sigma}^l\\[2mm]
&\quad+\|P\tilde{k}(\rho)\nabla^2\mathbf{u}\|_{\dot{B}_{2,\infty}^\sigma}^l+\|S\tilde{k}(\rho)\nabla^2\mathbf{u}\|_{\dot{B}_{2,\infty}^\sigma}^l.
\end{array}
\end{equation}
To deal with the terms on the right hand side of (\ref{4.14}), we use the same process as in (\ref{4.9})-(\ref{4.12}) to get
\begin{equation}\label{4.15}
\begin{array}{ll}
\|k'(0)P\nabla^2\mathbf{u}\|_{\dot{B}_{2,\infty}^\sigma}^l&\lesssim\|P^l\nabla^2\mathbf{u}\|_{\dot{B}_{2,\infty}^\sigma}^l
+\|P^h\nabla^2\mathbf{u}\|_{\dot{B}_{2,\infty}^\sigma}^l\\[2mm]
&\lesssim\|P^l\|_{\dot{B}_{2,\infty}^\sigma}\|\nabla^2\mathbf{u}\|_{\dot{B}_{p,1}^{\frac{n}{p}}}
+(\|\nabla^2\mathbf{u}\|_{\dot{B}_{p,1}^s}+\|\nabla^2\mathbf{u}^l\|_{\dot{B}_{2,1}^{\frac{n}{p}}})\|P^h\|_{\dot{B}_{p,1}^{-s}}\\[2mm]
&\lesssim\mathcal{E}_1(t)\|P^l\|_{\dot{B}_{2,\infty}^\sigma}+\mathcal{E}_\infty(t)\mathcal{E}_1(t),
\end{array}
\end{equation}
\begin{equation}\label{4.16}
\begin{array}{ll}
\|k'(0)S\nabla^2\mathbf{u}\|_{\dot{B}_{2,\infty}^\sigma}^l\!\!\!\!&\lesssim(\|S^l\|_{\dot{B}_{2,\infty}^\sigma}
\!+\!\|S^h\|_{\dot{B}_{2,1}^{\frac{n}{2}}})\|\mathbf{u}^l\|_{\dot{B}_{p,1}^{\frac{n}{p}+2}}
\!+\!\|S\|_{\dot{B}_{2,\infty}^{\frac{n}{p}-\frac{n}{2}+\sigma+1}}\!+\!\|\nabla^2\mathbf{u}^h\|_{\dot{B}_{p,1}^{\frac{n}{p}-1}}\\[3mm]
&\lesssim(\|S^l\|_{\dot{B}_{2,\infty}^\sigma}+\mathcal{E}_\infty(t))\mathcal{E}_1(t)+(\|S^l\|_{\dot{B}_{2,\infty}^\sigma}
+\|S^h\|_{\dot{B}_{2,1}^{\frac{n}{2}+1}})\|\mathbf{u}^h\|_{\dot{B}_{p,1}^{\frac{n}{p}+2}}\\[3mm]
&\lesssim\mathcal{E}_1(t)\|S^l\|_{\dot{B}_{2,\infty}^\sigma}+\mathcal{E}_\infty(t)\mathcal{E}_1(t),
\end{array}
\end{equation}
\begin{equation}\label{4.17}
\begin{array}{ll}
\|P\tilde{k}(\rho)\nabla^2\mathbf{u}\|_{\dot{B}_{2,\infty}^\sigma}^l
&\lesssim\|P^l\|_{\dot{B}_{2,\infty}^\sigma}\|\tilde{k}(\rho)\|_{\dot{B}_{p,1}^{\frac{n}{p}}}\|\nabla^2\mathbf{u}\|_{\dot{B}_{p,1}^{\frac{n}{p}}}
+(\|\tilde{k}(\rho)\|_{\dot{B}_{p,1}^{\frac{n}{p}}}\|\nabla^2\mathbf{u}\|_{\dot{B}_{p,1}^s}\\[3mm]
\!\!&\quad+\|\tilde{k}(\rho)\|_{\dot{B}_{p,1}^{\frac{n}{p}}}(\|\nabla^2\mathbf{u}^l\|_{\dot{B}_{2,1}^{\frac{n}{p}}}
+\|\nabla^2\mathbf{u}^h\|_{\dot{B}_{p,1}^{\frac{2n}{p}-\frac{n}{2}}}))\|P^h\|_{\dot{B}_{p,1}^{-s}}\\[3mm]
\!\!&\lesssim\mathcal{E}_\infty(t)\mathcal{E}_1(t)\|P^l\|_{\dot{B}_{2,\infty}^\sigma}+(\mathcal{E}_{\infty}(t))^2\mathcal{E}_1(t),
\end{array}
\end{equation}
\begin{equation}\label{4.18}
\begin{array}{ll}
\|S\tilde{k}(\rho)\nabla^2\mathbf{u}\|_{\dot{B}_{2,\infty}^\sigma}^l
&\lesssim\|S\|_{\dot{B}_{2,\infty}^\sigma}\|\tilde{k}(\rho)\|_{\dot{B}_{p,1}^{\frac{n}{p}}}\|\mathbf{u}^l\|_{\dot{B}_{p,1}^{\frac{n}{p}+2}}\\[4mm]
&\quad+(\|S\|_{\dot{B}_{p,1}^{\frac{n}{p}}}\|\tilde{k}(\rho)\|_{\dot{B}_{p,1}^s}
+\|\tilde{k}(\rho)\|_{\dot{B}_{p,1}^{\frac{n}{p}}}\|S\|_{\dot{B}_{2,1}^{\frac{n}{p}}})\|\nabla^2\mathbf{u}^h\|_{\dot{B}_{p,1}^{-s}}\\[4mm]
&\lesssim\mathcal{E}_\infty(t)\mathcal{E}_1(t)\|S^l\|_{\dot{B}_{2,\infty}^\sigma}+(\mathcal{E}_{\infty}(t))^2\mathcal{E}_1(t).
\end{array}
\end{equation}
Thus, we get
\begin{equation}\label{4.01}
\begin{array}{ll}
&\|k(\rho)(\mu\Delta \mathbf{u}+(\lambda+\mu)\nabla{\rm div}\mathbf{u}\|_{\dot{B}_{2,\infty}^\sigma}^l\\[2mm]
\lesssim&(1+\mathcal{E}_\infty(t)\mathcal{E}_1(t))\|(P,S)\|_{\dot{B}_{2,\infty}^\sigma}^l+(1+\mathcal{E}_{\infty}(t))\mathcal{E}_{\infty}(t)\mathcal{E}_1(t).
\end{array}
\end{equation}
Combining the above estimates, we infer from (\ref{4.2}) that
\begin{equation}\label{4.19}
\begin{array}{ll}
\|(P,\mathbf{u})\|_{\dot{B}_{2,\infty}^\sigma}^l&\displaystyle\lesssim\|(P_0,\mathbf{u}_0)\|_{\dot{B}_{2,\infty}^\sigma}^l+\int_0^t(1+\mathcal{E}_\infty(\tau))\mathcal{E}_\infty(\tau)\mathcal{E}_1(\tau)d\tau\\[2mm]
&\displaystyle\quad+\int_0^t(1+\mathcal{E}_\infty(\tau))\mathcal{E}_1(\tau)\|(P,\mathbf{u},S)\|_{\dot{B}_{2,\infty}^\sigma}^l d\tau.
\end{array}
\end{equation}
In the following, we deal with $\|S\|_{\dot{B}_{2,\infty}^\sigma}^l$. Applying $\dot\Delta_j$ to the third equation of (\ref{2.1}) gives
\begin{equation}\label{4.20}
\partial_t\dot\Delta_jS+\alpha_1\mathbf{u}\cdot\nabla\dot\Delta_jS+[\dot\Delta_j,\mathbf{u}\cdot\nabla]S=\frac{R}{\bar P}\dot\Delta_j((1-I(P))\Gamma(\alpha_1\mathbf{u})).
\end{equation}
Then,
\begin{equation}\label{4.21}
\begin{array}{ll}
\|S^l\|_{\dot B_{2,\infty}^\sigma}&\displaystyle\lesssim\|S_0^l\|_{\dot B_{2,\infty}^\sigma}+\int_0^t\|{\rm div}\mathbf{u}\|_{L^\infty}\|S^l\|_{\dot B_{2,\infty}^\sigma}d\tau\\[2mm]
&\displaystyle\quad+\int_0^t\|\nabla\mathbf{u}\|_{\dot B_{p,1}^{\frac{n}{p}}}\|S\|_{\dot B_{2,\infty}^\sigma}d\tau+\int_0^t\|(1-I(P))\Gamma(\alpha_1\mathbf{u})\|_{\dot B_{2,\infty}^\sigma}^ld\tau.
\end{array}
\end{equation}
With the aid of the embedding relation $\dot B_{p,1}^{\frac{n}{p}}\hookrightarrow L^\infty$, we can bound the first two terms on the right hand side of (\ref{4.21}) as
\begin{equation}\label{4.22}
\begin{array}{ll}
&\|{\rm div}\mathbf{u}\|_{L^\infty}\|S^l\|_{\dot B_{2,\infty}^\sigma}+\|\nabla \mathbf{u}\|_{\dot B_{p,1}^{\frac{n}{p}}}\|S\|_{\dot B_{2,\infty}^\sigma}\\[2mm]
\lesssim&\|\nabla \mathbf{u}\|_{\dot B_{p,1}^{\frac{n}{p}}}(\|S^l\|_{\dot B_{2,\infty}^\sigma}+\|S^h\|_{\dot B_{2,1}^{\frac{n}{2}}})\\[2mm]
\lesssim&\big(\|\mathbf{u}^l\|_{\dot B_{p,1}^{\frac{n}{p}+1}}+\|\mathbf{u}^h\|_{\dot B_{p,1}^{\frac{n}{p}+3}}\big)\big(\|S^l\|_{\dot B_{2,\infty}^\sigma}+\|S^h\|_{\dot B_{2,1}^{\frac{n}{2}}}\big)\\[2mm]
\lesssim&\mathcal{E}_1(t)\|S\|_{\dot B_{2,\infty}^\sigma}^l+\mathcal{E}_\infty(t)\mathcal{E}_1(t).
\end{array}
\end{equation}
For the last term in (\ref{4.21}), we get by a similar derivation of (\ref{4.5}) and (\ref{4.11}) that
\begin{equation}\label{4.23}
\begin{array}{ll}
\|(1-I(P))\Gamma(\alpha_1\mathbf{u})\|_{\dot B_{2,\infty}^\sigma}^l\!\!&\lesssim\|\Gamma(\alpha_1\mathbf{u})\|_{\dot B_{2,\infty}^\sigma}^l+\|I(P)\Gamma(\alpha_1\mathbf{u})\|_{\dot B_{2,\infty}^\sigma}^l\\[2mm]
\!\!&\lesssim\|\Gamma(\alpha_1\mathbf{u})\|_{\dot B_{2,\infty}^\sigma}^l+\|I(P)|\nabla \mathbf{u}|^2\|_{\dot B_{2,\infty}^\sigma}^l\\[2mm]
\!\!&\lesssim(1+\mathcal{E}_\infty(t))\mathcal{E}_1(t)\|\mathbf{u}\|_{\dot B_{2,\infty}^\sigma}^l+(1+\mathcal{E}_\infty(t))\mathcal{E}_\infty(t)\mathcal{E}_1(t).
\end{array}
\end{equation}
Inserting (\ref{4.22}) and (\ref{4.23}) into (\ref{4.21}) gives
\begin{equation}\label{4.24}
\begin{array}{ll}
\!\!\!\!\!\!\!\!\|S^l\|_{\dot B_{2,\infty}^\sigma}\!\lesssim\!\|S_0^l\|_{\dot B_{2,\infty}^\sigma}\!\!+\!\displaystyle\!\!\int_0^t(1\!+\!\mathcal{E}_\infty(\tau))\mathcal{E}_1(\tau)\|(\mathbf{u},S)\|_{\dot B_{2,\infty}^\sigma}^l\!\!d\tau\!+\!\!\int_0^t(1\!+\!\mathcal{E}_\infty(\tau))\mathcal{E}_\infty(\tau)\mathcal{E}_1(\tau)d\tau.
\end{array}
\end{equation}
Now, combining (\ref{4.19}) and (\ref{4.24}), one has
\begin{equation}\label{4.25}
\begin{array}{ll}
\displaystyle\|(P,\mathbf{u},S)\|_{\dot B_{2,\infty}^\sigma}^l\!\!&\displaystyle\lesssim\|(P_0,\mathbf{u}_0,S_0)\|_{\dot B_{2,\infty}^\sigma}^l+\int_0^t(1+\mathcal{E}_\infty(\tau))\mathcal{E}_\infty(\tau)\mathcal{E}_1(\tau)d\tau\\[2mm]
&\displaystyle\quad+\int_0^t(1+\mathcal{E}_\infty(\tau))\mathcal{E}_1(\tau)\|(P,\mathbf{u},S)\|_{\dot B_{2,\infty}^\sigma}^ld\tau.
\end{array}
\end{equation}
From Theorem \ref{Theorem 1.1}, we deduce that
\begin{equation}\label{4.26}
\begin{array}{ll}
\displaystyle\int_0^t(1+\mathcal{E}_\infty(\tau))\mathcal{E}_1(\tau)d\tau+\int_0^t(1+\mathcal{E}_\infty(\tau))\mathcal{E}_\infty(\tau)\mathcal{E}_1(\tau)d\tau\lesssim(1+\mathcal{X}_0)^3.
\end{array}
\end{equation}
Hence, by the Gronwall inequality, we obtain the desired uniform bound
\begin{equation}\label{4.27}
\begin{array}{ll}
\|(P,\mathbf{u},S)\|_{\dot B_{2,\infty}^\sigma}^l\leq C_0,\ for\ any\ \frac{n}{2}-\frac{2n}{p}\leq\sigma<\frac{n}{2}-1.
\end{array}
\end{equation}
Consequently, we complete the proof of Proposition \ref{proposition 4.1}.
\end{proof}
\end{proposition}

\subsection{Lyapunov-type differential inequality}
\quad\quad In this subsection, we develop the Lyapunov-type inequality in time for energy norms, which leads to the time-decay estimates. On the one hand, owing to $\frac{n}{2}-\frac{2n}{p}\leq\sigma<\frac{n}{2}-1$ and interpolation inequality, one has
\begin{equation}\label{4.32}
\|(P,\mathbf{u})\|_{\dot B_{2,1}^{\frac{n}{2}-1}}^l\leq C\big(\|(P,\mathbf{u})\|_{\dot B_{2,\infty}^\sigma}\big)^{\eta_1}\big(\|(P,\mathbf{u})\|_{\dot B_{2,1}^{\frac{n}{2}+1}}^l\big)^{1-\eta_1},\ \ \eta_1=\frac{4}{n-2\sigma+2}\in(0,1),
\end{equation}
which together with Proposition \ref{proposition 4.1} results that
\begin{equation}\label{4.33}
\|(P,\mathbf{u})\|_{\dot B_{2,1}^{\frac{n}{2}+1}}^l\geq c_0(\|(P,\mathbf{u})\|_{\dot B_{2,1}^{\frac{n}{2}-1}}^l)^{\frac{1}{1-\eta_1}}.
\end{equation}
On the other hand, due to $\mathcal{X}(t)<Cc_0$ and embedding relation in the high frequency, it is easy to see that
\begin{equation}\label{4.34}
\|\Lambda P\|_{\dot B_{p,1}^{\frac{n}{p}+1}}^h\geq C(\|\Lambda P\|_{\dot B_{p,1}^{\frac{n}{p}+1}}^h)^{\frac{1}{1-\eta_1}},\quad
\|\mathbf{u}\|_{\dot B_{p,1}^{\frac{n}{p}+3}}^h\geq C(\|\mathbf{u}\|_{\dot B_{p,1}^{\frac{n}{p}+1}}^h)^{\frac{1}{1-\eta_1}}.
\end{equation}
Thus, there exists a constant $\tilde{c}_0>0$ such that the following Lyapunov-type inequality holds
\begin{equation}\label{4.35}
\frac{d}{dt}(\|(P,\mathbf{u})^l\|_{\dot B_{2,1}^{\frac{n}{2}-1}}+\|(\Lambda P,\mathbf{u})\|_{\dot B_{p,1}^{\frac{n}{p}+1}}^h)+\tilde{c}_0(\|(P,\mathbf{u})^l\|_{\dot B_{2,1}^{\frac{n}{2}-1}}+\|(\Lambda P,\mathbf{u})\|_{\dot B_{p,1}^{\frac{n}{p}+1}}^h)^{1+\frac{4}{n-2\sigma-2}}\leq 0.
\end{equation}

\subsection{Decay estimate}
\quad\quad Solving the differential inequality (\ref{4.35}) directly, we get
\begin{equation}\label{4.36}
\|(P,\mathbf{u})^l\|_{\dot B_{2,1}^{\frac{n}{2}-1}}+\|(\Lambda P,\mathbf{u})\|_{\dot B_{p,1}^{\frac{n}{p}+1}}^h\leq C(1+t)^{-\frac{n-2\sigma-2}{4}}.
\end{equation}
For any $\sigma+n(\frac{1}{2}-\frac{1}{p})<\beta<\frac{n}{p}-1$, by the interpolation inequality we have
\begin{equation}\label{4.37}
\!\|(P,\mathbf{u})\|_{\dot B_{2,1}^{\beta+n(\frac{1}{2}-\frac{1}{p})}}^l\leq C(\|(P,\mathbf{u})\|_{\dot B_{2,\infty}^\sigma}^l)^{\eta_2}(\|(P,\mathbf{u})\|_{\dot B_{2,1}^{\frac{n}{2}-1}}^l)^{1-\eta_2},\quad\eta_2=\frac{\frac{n}{p}\!-\!1\!-\!\beta}{\frac{n}{2}\!-\!1\!-\!\sigma}\in(0,1),
\end{equation}
which combining with Proposition \ref{proposition 4.1} implies
\begin{equation}\label{4.38}
\|(P,\mathbf{u})\|_{\dot B_{2,1}^{\beta+n(\frac{1}{2}-\frac{1}{p})}}^l\leq C(1+t)^{-\frac{(\frac{n}{2}-\sigma-1)(1-\eta_2)}{2}}=C(1+t)^{-\frac{n}{2}(\frac{1}{2}-\frac{1}{p})-\frac{\beta-\sigma}{2}}.
\end{equation}
By virtue of $\sigma+n(\frac{1}{2}-\frac{1}{p})<\beta<\frac{n}{p}-1$, we find that
\begin{equation}\label{4.39}
\|(P^h,\mathbf{u}^h)\|_{\dot B_{p,1}^\beta}\leq C(\|P\|_{\dot B_{p,1}^{\frac{n}{p}+2}}^h+\|\mathbf{u}\|_{\dot B_{p,1}^{\frac{n}{p}+1}}^h)\leq C(1+t)^{-\frac{n-2\sigma-2}{4}},
\end{equation}
from which and (\ref{4.38}) yields
\begin{equation}\label{3.40}
\begin{array}{ll}
\!\!\!\!\!\|(P,\mathbf{u})\|_{\dot B_{p,1}^\beta}&\leq C(\|(P,\mathbf{u})\|_{\dot B_{2,1}^{\beta+n(\frac{1}{2}-\frac{1}{p})}}^l+\|(P,\mathbf{u})\|_{\dot B_{p,1}^\beta}^h)\\[2mm]
&\leq C(1+t)^{-\frac{n}{2}(\frac{1}{2}-\frac{1}{p})-\frac{\beta-\sigma}{2}}+C(1+t)^{-\frac{n-2\sigma-2}{4}}
\leq C(1+t)^{-\frac{n}{2}(\frac{1}{2}-\frac{1}{p})-\frac{\beta-\sigma}{2}}.
\end{array}
\end{equation}
Hence, thanks to the embedding relation $\dot B_{p,1}^0\hookrightarrow L^p$, one also has
\begin{equation}\label{3.41}
\begin{array}{ll}
\|\Lambda^\beta(P,\mathbf{u})\|_{L^p}\leq C(1+t)^{-\frac{n}{2}(\frac{1}{2}-\frac{1}{p})-\frac{\beta-\sigma}{2}}.
\end{array}
\end{equation}
The $L^p$-decay of $S_t$ can be derived by using the third equation of (\ref{1.5(0)}) and the decay result in (\ref{3.41}) directly. We omit the details for simplicity. This proves Theorem \ref{Theorem 1.2}.

\section {Appendix}

\quad\quad We will state several important lemmas and propositions on the homogeneous Besov space $\dot{B}^{s}_{p,1}$.
First, let $\mathcal {S}(\mathbb{R}^{n})$ be the Schwartz class of rapidly decreasing function. Given $f\in \mathcal
{S}(\mathbb{R}^{n})$, its Fourier transform $\mathcal {F}f=\widehat{f}$ is defined by $$\widehat{f}(\xi)=\int_{\mathbb{R}^{n}}e^{-ix\cdot\xi}f(x)dx.$$
Let $(\chi, \varphi)$ be a couple of smooth functions valued in $[0,1]$ such that $\chi$ is supported in the ball $\{\xi\in\mathbb{R}^{n}:  \ |\xi|\leq\frac{4}{3}\}$, $\varphi$ is supported in the shell $\{\xi\in \mathbb{R}^{n}: \ \frac{3}{4}\leq|\xi|\leq\frac{8}{3}\}$,  $\varphi(\xi):=\chi(\xi/2)-\chi(\xi)$
and
$$\chi(\xi)+\sum_{j\geq0}\varphi(2^{-j}\xi)=1~ \mathrm{for}\ \forall \ \xi \in\mathbb{ R}^{n},
~~~\sum_{j\in \mathbb{Z}}\varphi(2^{-j}\xi)=1 ~\mathrm{for} \ \forall \ \xi \in \mathbb{R}^{n}\setminus\{0\}.$$
For $f\in \mathcal{S}'$, the homogeneous frequency localization operators $\dot{\Delta}_j$ and $\dot{S}_j$ are defined by
\begin{equation*}
	\dot{\Delta}_{j}f\triangleq\varphi(2^{-j}D)f=\mathcal{F}^{-1}(\varphi(2^{-j}\xi)\mathcal{F}f)\quad {\rm and}\quad \dot{S}_{j}f\triangleq\chi(2^{-j}D)f=\mathcal{F}^{-1}(\chi(2^{-j}\xi)\mathcal{F}f).
\end{equation*}
We denote the space $\mathcal{S}'_{h}(\mathbb{R}^n)$ by the dual space of $\mathcal{S}'(\mathbb{R}^n)=\{f\in\mathcal{S}(\mathbb{R}^n):\,D^\alpha \hat{f}(0)=0\}$, which can also be identified by the quotient space of $\mathcal{S}'(\mathbb{R}^n)/{\mathbb{P}}$ with the polynomial space ${\mathbb{P}}$. The formal equality $$ f=\sum_{j\in\mathbb{Z}}\dot{\Delta}_jf $$ holds true for $f\in\mathcal{S}'_{h}(\mathbb{R}^n)$ and is called the homogeneous Littlewood-Paley decomposition, and then we have the fact that
\begin{equation}
	\dot{S}_jf=\sum_{q\le j-1}\dot{\Delta}_qf. \nonumber
\end{equation}
One easily verifies that with our choice of $\varphi$,
\begin{equation*}
	\dot{\Delta}_j\dot{\Delta}_qf\equiv0\quad \textrm{if}\quad|j-q|\ge
	2\quad \textrm{and} \quad
	\dot{\Delta}_j(\dot{S}_{q-1}f\dot{\Delta}_q f)\equiv0\quad \hbox{if}
	\quad |j-q|\ge 5.
\end{equation*}

\begin{definition}(Homogeneous Besov space)\label{A.1}
	For $s\in \mathbb{R}$ and $1\le p,r\le \infty$, the homogeneous Besov space $\dot{B}^{s}_{p,1}$ is defined by
	\begin{equation}\label{a.1}
		\dot{B}^s_{p,r}\triangleq\left\{f\in \mathcal{S}_h':||f||_{\dot{B}^s_{p,r}}<+\infty\right\},
	\end{equation}
	where
	\begin{equation}\label{a.2}
		||f||_{\dot{B}^s_{p,r}}\triangleq||2^{js}||\dot{\Delta}_jf||_{L^p}||_{\mathit{l}^r(\mathbb{Z})}.
	\end{equation}
\end{definition}

\begin{definition}(Chemin-Lerner spaces)\label{A.2}
	Let $T>0$, $s\in\mathbb{R}$, $1<r,p,q\le\infty$. The space $\widetilde{L}^q_{T}(\dot{B}^s_{p,r})$ is defined by
	\begin{equation}\label{a.3}
		\widetilde{L}^q_{T}(\dot{B}^s_{p,r})\triangleq\left\{f\in L^q(0,T;\mathcal{S}'_h):||f||_{\widetilde{L}^q_{T}(\dot{B}^s_{p,r})}<+\infty\right\},
	\end{equation}
	where
	\begin{equation}\label{a.4}
		||f||_{\widetilde{L}^q_{T}(\dot{B}^s_{p,r})}\triangleq||2^{js}||\dot{\Delta}_jf||_{L^q(0,T;L^p)}||_{\mathit{l}^r(\mathbb{Z})}.
	\end{equation}
\end{definition}
\begin{remark}\label{A.3}
	It holds that
	\begin{equation*}
		||f||_{\widetilde{L}^q_{T}(\dot{B}^s_{p,r})}\le||f||_{L^q_{T}(\dot{B}^s_{p,r})}\quad {\rm if}\quad r\ge q;\quad||f||_{\widetilde{L}^q_{T}(\dot{B}^s_{p,r})}\ge||f||_{L^q_{T}(\dot{B}^s_{p,r})}\quad {\rm if}\quad r\le q.
	\end{equation*}
\end{remark}
Restricting the above norms (4.2), (4.3) to the low or high frequencies parts of distributions will be crucial in our approach. For example, let us fix some integer $j_0$ and set
\begin{equation*}
	||f||^l_{\dot{B}^{s}_{p,1}}\triangleq\sum_{j\le j_0}2^{js}||\dot{\Delta}_jf||_{L^p},\quad ||f||^h_{\dot{B}^{s}_{p,1}}\triangleq\sum_{j\ge j_0-1}2^{js}||\dot{\Delta}_jf||_{L^p};
\end{equation*}
\begin{equation*}
	||f||^l_{\widetilde{L}^\infty_{T}(\dot{B}^s_{p,1})}\triangleq\sum_{j\le j_0}2^{js}||\dot{\Delta}_jf||_{L^\infty_T(L^p)},\quad ||f||^h_{\widetilde{L}^\infty_{T}(\dot{B}^s_{p,1})}\triangleq\sum_{j\ge j_0-1}2^{js}||\dot{\Delta}_jf||_{L^\infty_T(L^p)}.
\end{equation*}
\begin{lemma}(Bernstein inequalities)\label{A.4}
	Let $\mathscr{B}$ be a ball and $\mathscr{C}$ be a ring of $\mathbb{R}^n$. For $\lambda>0$, integer $k\ge0$, $1\le p\le q\le \infty$ and a smooth homogeneous function $\sigma$ in $\mathbb{R}^n\backslash\{0\}$ of degree $m$, then there holds
	\begin{equation*}
		||\nabla^kf||_{L^q}\le C^{k+1}\lambda^{k+n(\frac{1}{p}-\frac{1}{q})}||f||_{L^p},\quad {\rm whenever\ supp}\widehat{f}\subset\lambda\mathscr{B},
	\end{equation*}
	\begin{equation*}
		C^{-k-1}\lambda^k||f||_{L^q}\le||\nabla^kf||_{L^p}\le C^{k+1}\lambda^k||f||_{L^p},\quad {\rm whenever\ supp}\widehat{f}\subset\lambda\mathscr{C},
	\end{equation*}
	\begin{equation*}
		||\sigma(\nabla)f||_{L^q}\le C_{\sigma,m}\lambda^{m+n(\frac{1}{p}-\frac{1}{q})}||f||_{L^p},\quad {\rm whenever\ supp}\widehat{f}\subset\lambda\mathscr{C}.
	\end{equation*}
\end{lemma}

\begin{proposition}\cite{Bahouri-2011}(Embedding for Besov space on $\mathbb{R}^n$)\label{A.5}
	\begin{itemize}
		\item For any $p\in[1,\infty]$, we have the continuous embedding $\dot{B}^{0}_{p,1}\hookrightarrow L^p\hookrightarrow\dot{B}^{0}_{p,\infty}$.
		\item If $s\in\mathbb{R}$, $1\le p_1\le p_2\le \infty$, and $1\le r_1\le r_2\le \infty$ then  $\dot{B}^{s}_{p_1,r_1}\hookrightarrow\dot{B}^{s-n(\frac{1}{p_1}-\frac{1}{p_2})}_{p_2,r_2}$.
		\item The space $\dot{B}^{\frac{n}{p}}_{p,1}$ is continuously embedded in the set of bounded continuous function (going to zero at infinity if, additionally, $p<\infty$).
	\end{itemize}
\end{proposition}
\begin{proposition}\cite{duan3}\label{A.6}
	If supp$\mathcal{F}f\subset\left\{\xi\in\mathbb{R}^n:R_1\lambda\le|\xi|\le R_2\lambda\right\}$, then there exists $C$ depending only on $d$, $R_1$, $R_2$ so that for all $1<p<\infty$,
	\begin{equation}\label{a.5}
		C\lambda^2(\frac{p-1}{p})\int_{\mathbb{R}^n}|f|^pdx\le(p-1)\int_{\mathbb{R}^n}|\nabla f|^2|f|^{p-2}dx=-\int_{\mathbb{R}^n}\Delta f|f|^{p-2}fdx.
	\end{equation}
\end{proposition}
\begin{proposition}\cite{Bahouri-2011}(Interpolation inequality)\label{interpolation}
	Let $1\le p,r,r_1,r_2\le\infty$, if $f\in\dot{B}^{s_1}_{p,r_1}\cap\dot{B}^{s_2}_{p,r_2}$ and $s_1\neq s_2$, then $f\in\dot{B}^{\theta s_1+(1-\theta)s_2}_{p,r}$ for all $\theta\in(0,1)$ and
	\begin{equation}\label{a.6}
		||f||_{\dot{B}^{\theta s_1+(1-\theta)s_2}_{p,r}}\le||f||^\theta_{\dot{B}^{s_1}_{p,r_1}}||f||^{1-\theta}_{\dot{B}^{s_2}_{p,r_2}}
	\end{equation}
	with $\frac{1}{r}=\frac{\theta}{r_1}+\frac{1-\theta}{r_2}$.
\end{proposition}

\begin{proposition}\cite{Bahouri-2011,danchin5}\label{A.7}
	Let $s>0$, $1\leq p$, $r\leq \infty$, then $\dot{B}^{s}_{p,r}\cap L^\infty$ is an algerbra and
	\begin{equation}\label{a.7}
		\begin{array}{rl}
			||fg||_{\dot{B}^{s}_{p,r}} \lesssim ||f||_{L^\infty}||g||_{\dot{B}^{s}_{p,r}} + ||g||_{L^\infty}||f||_{\dot{B}^{s}_{p,r}}.
		\end{array}
	\end{equation}
	
	Let $s_1+s_2>0$, $s_1 \le \frac{d}{p_1}$, $s_2 \le \frac{d}{p_2}$, $s_1\ge s_2$, $\frac{1}{p_1}+\frac{1}{p_2}\le 1$. Then it holds that
	\begin{equation}\label{a.8}
		||fg||_{\dot{B}^{s_2}_{q,1}} \lesssim ||f||_{\dot{B}^{s_1}_{p_1,1}}||g||_{\dot{B}^{s_2}_{p_2,1}},
	\end{equation}
	where $\frac{1}{q} = \frac{1}{p_1} + \frac{1}{p_2} - \frac{s_1}{d}$.
\end{proposition}
\begin{proposition}\cite{xin3}\label{A.8}
	Let the real numbers $s_1,\ s_2,\ p_1$ and $p_2$ be such that
	$$s_1+s_2\geq0,\ s_1\leq\frac{n}{p_1},\ s_2<\min\left(\frac{n}{p_1},\frac{n}{p_2}\right)\ and\ \frac{1}{p_1}+\frac{1}{p_2}\leq1.$$
	Then it holds that
	\begin{equation}\label{a.9}
		\|fg\|_{\dot B_{p_2,\infty}^{s_1+s_2-\frac{n}{p_1}}}\lesssim\|f\|_{\dot B_{p_1,1}^{s_1}}\|g\|_{\dot B_{p_2,\infty}^{s_2}}.
	\end{equation}
\end{proposition}

\begin{corollary}\label{A.9}
	Let the real numbers $1-\frac{n}{2}<\sigma_1\leq\sigma_0$ and $p$ satisfy (\ref{1.6}). The following two inequalities hold true:
	\begin{equation}\label{a.10}
		\|fg\|_{\dot B_{2,\infty}^{-\sigma_1}}\lesssim\|f\|_{\dot B_{p,1}^{\frac{n}{p}}}\|g\|_{\dot B_{2,\infty}^{-\sigma_1}},
	\end{equation}
	as well as
	\begin{equation}\label{a.11}
		\|fg\|_{\dot B_{2,\infty}^{\frac{n}{p}-\frac{n}{2}-\sigma_1}}\lesssim\|f\|_{\dot B_{p,1}^{\frac{n}{p}-1}}\|g\|_{\dot B_{2,\infty}^{\frac{n}{p}-\frac{n}{2}-\sigma_1+1}}.
	\end{equation}
\end{corollary}

\begin{proposition}\cite{xu2023}\label{A.10}
	Let $j_0\in\mathbb{Z}$, and denote $z^l\triangleq\dot S_{j_0}z,\ z^h\triangleq z-z^l$ and, for any $s\in\mathbb{R}$,
	$$\|z\|_{\dot B_{2,\infty}^{s}}^l\triangleq \sup_{j\leq j_0}2^{js}\|\dot\Delta_jz\|_{L^2}.$$
	There exists a universal integer $N_0$ such that for any $2\leq p\leq 4$ and $s>0$, we have
	\begin{equation}\label{a.12}
		\|fg^h\|_{\dot B_{2,\infty}^{-\sigma_0}}^l\leq C(\|f\|_{\dot B_{p,1}^{s}}+\|\dot S_{k_0+N_0}f\|_{L^{p^\ast}})\|g^h\|_{\dot B_{p,\infty}^{-s}},
	\end{equation}
	\begin{equation}\label{a.13}
		\|f^h g\|_{\dot B_{2,\infty}^{-\sigma_0}}^l\leq C(\|f^h\|_{\dot B_{p,1}^{s}}+\|\dot S_{k_0+N_0}f^h\|_{L^{p^\ast}})\|g\|_{\dot B_{p,\infty}^{-s}}.
	\end{equation}
	
	Additionally, for exponents $s>0$, $1 \le p_1,p_2,q \le \infty$ satisfying
	\begin{equation*}
		\frac{d}{p_1}+\frac{d}{p_2}-d\le s\le \min (\frac{d}{p_1},\frac{d}{p_2})\ \ and\ \ \frac{1}{q} = \frac{1}{p_1}+\frac{1}{p_2}-\frac{s}{d}.
	\end{equation*}
	Then it holds that
	\begin{equation}\label{a.14}
		||fg||_{\dot{B}^{-s}_{q,\infty}} \lesssim ||f||_{\dot{B}^{s}_{p_1,1}}||g||_{\dot{B}^{-s}_{p_2, \infty}}.
	\end{equation}
\end{proposition}

\begin{proposition}\cite{xu2023}\label{A.11} The Bony decomposition satisfies that
	\begin{equation}\label{a.15}
		||T_ab||_{\dot{B}^{s-1+ \frac{d}{2} - \frac{d}{p}}_{2,1}} \lesssim ||a||_{_{\dot{B}^{\frac{d}{p}-1}_{p,1}}}||b||_{\dot{B}^{s}_{p,1}},\ \ if\ d \ge 2\ and\ 1 \le p \le \min(4,\frac{2d}{d-2}),
	\end{equation}
	\begin{equation}\label{a.16}
		||R(a,b)||_{\dot{B}^{s-1+ \frac{d}{2} - \frac{d}{p}}_{2,1}} \lesssim ||a||_{_{\dot{B}^{\frac{d}{p}-1}_{p,1}}}||b||_{\dot{B}^{s}_{p,1}},\ \ if\ s>1-\min(\frac{d}{p}+\frac{d}{p'})\ and\ 1\le p\le 4.
	\end{equation}
\end{proposition}

\begin{lemma}\cite{danchin5}\label{A.12}
	Let $n\geq2$, $1\le p, q\le \infty$, $v\in \dot{B}^{s}_{q,1}(\mathbb{R}^n)$ and $\nabla u \in \dot{B}^{\frac{n}{p}}_{p,1}(\mathbb{R}^n)$.\\
	Asumme that
	\begin{equation*}
		-n \min(\frac{1}{p}, 1-\frac{1}{q})<s\le 1+n\min(\frac{1}{p}, \frac{1}{q}).
	\end{equation*}
	Then it holds the commutator estimate
	\begin{equation}\label{a.17}
		||[\dot{\Delta}_j, u\cdot\nabla]v||_{L^p} \lesssim d_j2^{-js}||\nabla u||_{\dot{B}^{\frac{n}{p}}_{p,1}}||v||_{\dot{B}^{s}_{q,1}}.
	\end{equation}
	In the limit case $s=-n\min(\frac{1}{p}, 1-\frac{1}{q})$, we have
	\begin{equation}\label{a.18}
		\sup2^{js}||[\dot{\Delta}_j, u\cdot\nabla]||_{L^p}\lesssim ||\nabla u||_{\dot{B}^{\frac{n}{p}}_{p,1}}||v||_{\dot{B}^{s}_{q,\infty}}.
	\end{equation}
\end{lemma}

\begin{proposition}\cite{Bahouri-2011}\label{A.13}
	Let $F: \mathbb{R} \mapsto \mathbb{R}$ be a smooth function with $F(0)=0$, $1\le p, r\le \infty$ and $s>0$. Then $F: \dot{B}^{s}_{p,r}(\mathbb{R}^n)\cap(L^\infty(\mathbb{R}^n))$ and
	\begin{equation}\label{a.19}
		||F(u)||_{\dot{B}^{s}_{p,r}}\le C||u||_{\dot{B}^{s}_{p,r}}
	\end{equation}
	with $C$ a constant depending only on $|||u||_{L^\infty}$, $s$, $p$, $n$ and derivatives of $F$.
	
	If $s>-\min(\frac{n}{p}, \frac{n}{p'})$, then $F: \dot{B}^{s}_{p,r}(\mathbb{R}^n)\cap\dot{B}^{\frac{n}{p}}_{p,1}(\mathbb{R}^n) \mapsto \dot{B}^{s}_{p,r}(\mathbb{R}^n)\cap\dot{B}^{\frac{n}{p}}_{p,1}(\mathbb{R}^n)$, and
	\begin{equation}\label{a.20}
		||F(u)||_{\dot{B}^{s}_{p,r}}\le C(1+||u||_{\dot{B}^{\frac{n}{p}}_{p,1}})||u||_{\dot{B}^{s}_{p,r}}.
	\end{equation}
\end{proposition}

		\section*{Acknowledgments}
		
		\bigbreak
		
		{\bf Funding}: The research was supported by National Natural Science Foundation of China (11971100) and Natural Science Foundation of Shanghai (22ZR1402300).\\
		{\bf Conflict of Interest}: The authors declare that they have no conflict of interest.

\bibliographystyle{plain}

	\end{document}